\documentclass[final,12pt]{colt2019} 

\title[PAC-Bayesian Bound for DGPs]{PAC-Bayesian Bounds for Deep Gaussian Processes}
\usepackage{times}
\usepackage{float}
\usepackage{microtype}
\usepackage{booktabs} 
\usepackage[utf8]{inputenc}
\usepackage{upgreek}
\usepackage{mathtools}
\usepackage{url}
\usepackage{textgreek}
\usepackage{multirow}
\usepackage{scrextend}
\usepackage{algorithm}
\usepackage{algpseudocode}
\usepackage{lineno}
\usepackage[overload]{empheq}
\usepackage{arydshln}
\DeclareMathAlphabet{\mathpzc}{OT1}{pzc}{m}{it}
\DeclareFontFamily{OT1}{pzc}{}
\DeclareFontShape{OT1}{pzc}{m}{it}{<-> s * [1.200] pzcmi7t}{}
\DeclareMathAlphabet{\mathpzc}{OT1}{pzc}{m}{it}
\newtheorem{defi}{Definition}

\newcommand{\subalign}[2][c]{%
  \if#1c\vcenter\else\vtop\fi{%
    \ialign{\hfil$\m\scriptstyle##$&$\m\scriptstyle{}##$\crcr
      #2\crcr
    }%
  }
}
\makeatother
\newcommand\numberthis{\addtocounter{equation}{1}\tag{\theequation}}
\DeclareMathAlphabet\mathbfcal{OMS}{cmsy}{b}{n}
\DeclareMathAlphabet{\mathbbmsl}{U}{bbm}{m}{sl}

\coltauthor{\Name{Roman Föll},  \Email{foell@mathematik.uni-stuttgart.de}\\
\addr Institute of Applied Analysis and Numerical Simulation\\
University of Stuttgart\\
Stuttgart, 70569, Germany\\
\AND
\Name{Ingo Steinwart}, \Email{steinwart@mathematik.uni-stuttgart.de}\\
\addr Institute of Stochastics and Applications\\
University of Stuttgart\\
Stuttgart, 70569, Germany\\
}

\begin{document}

\maketitle

\begin{abstract}%
Variational approximation techniques and inference for stochastic models in machine learning has gained much attention the last years. Especially in the case of Gaussian Processes (GP) and their deep versions, Deep Gaussian Processes (DGPs), these viewpoints improved state of the art work. In this paper we introduce Probably Approximately Correct (PAC)-Bayesian risk bounds for DGPs making use of variational approximations. We show that the minimization of PAC-Bayesian generalization risk bounds maximizes the variational lower bounds belonging to the specific DGP model. We generalize the loss function property of the log likelihood loss function in the context of PAC-Bayesian risk bounds to the quadratic-form-Gaussian case. Consistency results are given and an oracle-type inequality gives insights in the convergence between the raw model (predictor without variational approximation) and our variational models (predictor for the variational approximation). Furthermore, we give extensions of our main theorems for specific assumptions and parameter cases. Moreover, we show experimentally the evolution of the consistency results for two Deep Recurrent Gaussian Processes (DRGP) modeling time-series, namely the recurrent Gaussian Process (RGP) and the DRGP with Variational Sparse Spectrum approximation, namely DRGP-(V)SS. \end{abstract}

\begin{keywords}%
  PAC-Bayesian theory, Deep Gaussian Process models, Variational approximations, Consistency, Recurrent models
\end{keywords}

\section{Introduction}
\label{sec:introduction}

The Bayesian viewpoint for probabilistic inference is very popular in the statistics and the machine learning communities~\cite{rasmussen2006gaussian,neal2012bayesian,kingma2013auto}. Its flexibility and simple framework are important factors for their success. On the one hand, regarding many applications, Bayesian approaches represent state of the art benchmark methods~\cite{al2016learning, salimbeni2017doubly}, on the other hand, the PAC-Bayesian approach is a powerful way to derive risk bounds for probabilistic models generated from Bayesian modeling and inference. It originates from~\cite{Shawe-Taylor, catoni2004statistical, catoni2007pac, mcallester1999some, mcallester1999pac}. In this paper, the focus is on PAC-Bayesian investigations, that use variational approximations and inference instead of Bayesian inference. Variational approximation and inference is a promising approach in the last few years in many fields. We will focus on recent work of~\cite{alquier2015properties, germain2016pac, NIPS20177100} to derive PAC-Bayesian results for a class of variational stochastic models. More precisely, we derive PAC-Bayesian statements for the DGP models of~\cite{damianou2015deep, mattos2015recurrent, cutajar2016practical, foell2019deep}, which use variational approximations and inference instead of Bayesian inference. Unlike simple GPs, DGPs have proven to be capable of capturing non-stationarity and heteroscedasticity inherent to many modeling problems and applications in practice. The theoretical aspect, regarding generalization properties of DGPs, is until now rather less understood. We can use PAC-Bayesian theory to provide precise answers by a guaranteed upper bound on the generalization error in an unspecified data-distribution setting with high probability. Therefore, using these statements in practice, we are able to design stochastic models, here DGPs, that have good generalization properties with high probability. In our setting for the supervised learning case of regression, we assume, that the output-data is coming from a multi-variate distribution given the input-data. Moreover, we assume a fixed design scenario, where the input-variables are set by an experimenter. It occurs in many practical applications like in the controlling or the prediction/simulation case and is the standard in regression~\cite{deisenroth2013gaussian,al2016learning}.

\section{Related Work to theoretical analysis of DGPs}
\label{sec:RelatedWork}

Early work on studying GPs in the PAC-Bayesian approach goes back to~\cite{seeger2002pac, seeger2003bayesian} for classifiers. Further developments have been made by~\cite{van2008rates, vaart2011information} who investigated the convergence rate of GP estimators regarding geometric relations between the true function and the \textit{Reproducing kernel Hilbert space} (RKHS) corresponding to the GP prior. Based on this,~\cite{suzuki2012pac} developed PAC-Bayesian oracle inequalities for GP regression and Multiple Kernel Additive Models with convergence results, where they could improve some of the previous results of~\cite{van2008rates, vaart2011information}. Regarding \textit{deep neural networks} (DNN),~\cite{duvenaud2014avoiding} studied DGPs, a type of infinitely-wide DNN, see~\cite{lee2017deep}, and deep kernels, as well as their pathologies and how these pathologies could be alleviated. Recently,~\cite{dunlop2017deep} developed a unifying perspective on hierarchical GPs, leading to a wide class of DGPs. Exploiting the fact, that this common framework has a Markovian structure, they interpret the depth of the process in terms of the ergodicity or non-ergodicity of this process. Their analysis is based solely on the DGP for unobserved data, and not the conditioned process in the inference problem with observed data. Our derived PAC-Bayesian statements for the empirical bound case and the oracle-type case are valid for both kind of scenarios. \\
To our knowledge, we are the first to derive explicit PAC-Bayesian statements for the DGP models of \cite{damianou2015deep, mattos2015recurrent, cutajar2016practical, foell2019deep}. As mentioned in the introduction,~\cite{alquier2015properties, germain2016pac, NIPS20177100} present a solid basis within the Bayesian and variational framework to derive our new results.

\section{PAC-Bayesian Theory: Notation and Definitions}
\label{sec:NotationandDefinitions}

In the following, we use the notation $\mathpzc{f}$, $\mathpzc{y}$ for stochastic processes, $f_\mathbf{x}$, $\boldsymbol{y}$, (italic) for random variables, $\mathrm{f}(\mathbf{x})$, $\mathbf{y}$ (upright) for realizations and data. We assume, that we are given a \textit{bounded} set of input-states
\begin{align*}
\mathbf{x}_\mathsf{1},\dots,\mathbf{x}_\mathsf{K}\in\mathbb R^Q,\quad \mathbf{X}\stackrel{\mathrm{def}}=[\mathbf{x}_1,\dots,\mathbf{x}_\mathsf{K}]^T\in\mathbfcal{X}\stackrel{\mathrm{def}}=\mathbb{R}^{\mathsf{K}\times Q}, 
\end{align*}
where $\mathsf{K}\geq 1$ and $Q\in\mathbb N$. Moreover, we assume, that conditional on $\mathbf{X}$, there is some unknown multivariate data-distribution $P_{\mathsf{K}}$ on $\mathbb R^{\mathsf{K}}$ which generates $N\in\mathbb N$ observations
\begin{flalign}
\mathbf{y}^\mathsf{1}_{\mathbf{X}},\dots,\mathbf{y}^N_{\mathbf{X}}\in\mathbfcal{Y}\stackrel{\mathrm{def}}= \mathbb R^{\mathsf{K}},\quad \mathbb{D}\stackrel{\mathrm{def}}=\{\mathbf{y}^i_{\mathbf{X}}\}_{i=1}^N,\numberthis\label{datapac}
\end{flalign}
\textit{independently, identically distributed} (iid), and we write $\mathbf{y}^i\stackrel{\mathrm{def}}=\mathbf{y}^i_{\mathbf{X}}\stackrel{\mathrm{def}}=[\mathrm{y}^{i}_{\mathbf{x}_1},\dots,\mathrm{y}^{i}_{\mathbf{x}_\mathsf{K}}]^T\in\mathbb{R}^{\mathsf{K}}$ as a shorthand. The multi-variate data-distribution $P_{\mathsf{K}}$ might be given as the marginal distribution of a stochastic process $\mathpzc{y}$, which we do not know. For consistency reasons, modeling with a specific DGP, this is reasonable to assume. Remember, a DGP is an stochastic process build by stacking GPs. A unbounded stochastic process $\mathpzc{y}\stackrel{\mathrm{def}}=[y_{\mathbf{x}}]_{\mathbf{x}\in\mathbb R^Q}$ is a GP if and only if any finite collection of random variables $\boldsymbol{y}_{\mathbf{X}}\stackrel{\mathrm{def}}=[y_{\mathbf{x}_1},\dots,y_{\mathbf{x}_\mathsf{K}}]^T$ forms a multivariate Gaussian random vector~\citep{rasmussen2006gaussian}. Because a GP is unbounded, we are not allowed to e.g. restrict the components $y_{\mathbf{x}_i}$ to be bounded. Moreover, the model parameters $\uptheta$, applied for the modeling task with DGPs, are used to prove our results. We want to emphasize here, that we do assume the standard case, where the sample-vectors $\mathbf{y}^\mathsf{1},\dots,\mathbf{y}^N$ are always iid, but the entries $\mathrm{y}^{i}_{\mathbf{x}_1},\dots,\mathrm{y}^{i}_{\mathbf{x}_\mathsf{K}}$ of the real-vector $\mathbf{y}^i_{\mathbf{X}}$ might be dependent. So note, $P_{\mathsf{K}}$ is in general not a product distribution. 
We also write $\boldsymbol{y}_{\mathbf{X},\uptheta}$ for the multivariate random vector depending on our model parameters. If the context is clear, we also write $\boldsymbol{y}$.\\ Furthermore, we have an independent output data-set $\mathbf{Y}\in\mathbb R^{\mathsf{K}\times \bar{N}}$ for the training task with DGP models, where $\bar{N}\in\mathbb N$ is independent from the sampling amount $N$ and $\bar{N}\ll N$. This output data-set $\mathbf{Y}$ is assumed to be observed on the same states $\mathbf{X}$ as the observations $\mathbb{D}$, but the difference is, that $\mathbf{Y}$ is assumed to be constant (we condition on these) to derive the posterior distribution of the model. This setting and the definitions make sense for both applications, the simple static regression case and the dynamic modeling case (modeling time series), as $\mathbf{x}$ represents always some state, e.g. state of time or some physical state. The choice of the amount of states $\mathbf{x}$ given by $\mathsf{K}$ should be seen as an experimenters choice for his specific modeling task. These predefined states $\mathbf{X}$ are observed $N$ times in our PAC-Bayesian framework and our goal is, that the generalization error tends to zero, as $N$ tends to infinity. Additionally, we want to emphasize that many experimenters in practice often choose arbitrary states and collect just single measurements on these. The choice of modeling with a DGP involves measurement errors for outputs $\mathrm{y}^{i}_{\mathbf{x}_\mathsf{k}}$ at specific states $\mathbf{x}_\mathsf{k}$, which should be therefore measured several times. Nevertheless, this hints to a proper data collection in theory.\\
Furthermore, we have $\mathfrak{f}_{\uptheta}:\mathbfcal{X}\times\Theta\to\mathbfcal{Y}$ as our predictor and where $\mathcal{F}_\Theta\stackrel{\mathrm{def}}=\{\mathfrak{f}_{\uptheta}:\mathbfcal{X}\times\Theta \to\mathbfcal{Y},\uptheta\in\Theta\}$ is the set of all these predictors, whereby we have $\uptheta\in\Theta\subset\mathbb R^p$ and $\mathbb R^p$, $p\in\mathbb N$, is a $p$-dimensional parameter space. The size of $p$ depends on our chosen model and the specific parameter families. The set $\Theta$ represents the restriction to the ranges of the specific parameter families for our models. Do not confuse $\mathfrak{f}_{\uptheta}$ with the marginal, noise free random vector $\textit{\textbf{f}}_{\mathbf{X},\uptheta}$ of the noise free GP $\mathpzc{f}$. In a regression context we assume $\mathbf{y}^{i} = \mathrm{f}(\mathbf{X})+\epsilon_i^{\mathbf{y}}$, where the respective function values $\mathbf{f}=\mathrm{f}(\mathbf{X})=[\mathrm{f}(\mathbf{x}_1),\dots,\mathrm{f}(\mathbf{x}_n)]^T$ are not the predictor for a GP (for a GP, we use the mean-function as predictor). Moreover, we consider an unbounded loss function $\ell:\mathcal F_{\Theta}\times\mathbfcal{Y}\to\mathbb R$, and denote $\mathcal{G}_{\Theta}$ as the sets of all probability distributions on the model parameters. To link the PAC-Bayesian theory and the variational framework to the regression context, we have to define the empirical risk and the generalization error. We write the empirical risk and the generalization error as 
\begin{align*} 
\mathcal{L}_{\mathbb{D}}^{\ell}(\mathfrak{f}_{\uptheta})\stackrel{\mathrm{def}}=\frac{1}{N}\sum\nolimits_{i=1}^{N} \ell(\mathfrak{f}_{\uptheta},\mathbf{y}^i),\quad\text{and}\quad\mathcal{L}_{P_\mathsf{K}}^{\ell}(\mathfrak{f}_{\uptheta})\stackrel{\mathrm{def}}=\underset{\boldsymbol{y}\sim P_\mathsf{K}}{\mathbf{E}}\left[\ell(\mathfrak{f}_{\uptheta},{\mathbf{y}})\right].
\end{align*}
Furthermore, we define for $\lambda>0$ the Gibbs posterior $g_{\lambda}$, a density
\begin{equation*} 
g_{\lambda}(\uptheta)\stackrel{\mathrm{def}}=\frac{e^{-\lambda\mathcal{L}_{\mathbb{D}}^{\ell}(\mathfrak{f}_{\uptheta})}\pi(\hat{\uptheta})}{\int e^{-\lambda\mathcal{L}_{\mathbb{D}}^{\ell}(\mathfrak{f}_{\uptheta})}\pi(\hat{\uptheta})d\hat{\uptheta}},\quad \hat\uptheta\in\hat{\Theta}\subset\Theta,
\end{equation*}
where $\pi$ is a prior density over $\hat\uptheta$. We further have $\mathcal{G}_{\hat{\Theta}}$, which is the corresponding set of all probability distributions on these model parameters. Moreover $\check{\Theta}\subset\check{\boldsymbol\Theta}$ is the subset of all parameters without a prior assumption and ${\uptheta}_\mathrm{m}$, ${\uptheta}_\mathrm{v}$ are the variational parameters for $\hat\uptheta$ (mean, variance). We will define for every GP from $l=1,\dots,L+1$ in the DGP separate $\uptheta^{(l)}$ in the Appendix~\ref{pacbayesian}, Equation~(\ref{predictor1})-~(\ref{predictor3}) and $\uptheta$ will be the stacked version of them. Let as summarize the definitions of our parameter sets as
\begin{align*}
\uptheta=(\hat\uptheta,\check\uptheta)^T=(\hat\uptheta^{(1)},\check\uptheta^{(1)},\dots,\hat\uptheta^{(L+1)},\check\uptheta^{(L+1)})^T\in\Theta\subset\mathbb R^p.
\end{align*}
As we will later see, $\mathcal{G}_{\hat{\Theta}}$ also represents the space of our variational distributions. In the case of choosing the loss as the negative log likelihood $\ell=\ell_{\text{nll}} \stackrel{\mathrm{def}}= -\log(p(\mathbf{y}|\uptheta,\mathbf{X}))$ and $\lambda = N$, the Gibbs posterior coincides with the Bayesian posterior, see~\cite{germain2016pac}. We assume the best possible variational approximation on $\hat{\Theta}$ exists and is 
\begin{equation*} 
Q_{\lambda}\stackrel{\mathrm{def}}=\text{arg}\min\limits_{Q_{\text{\tiny{PAC}}}\in\mathcal{G}_{\hat{\Theta}}}\mathbf{KL}(Q_{\text{\tiny{PAC}}}||g_{\lambda}),
\end{equation*}
where $\mathbf{KL}$ denotes the Kullback-Leibler divergence. The expression $Q_{N}$ is then the best variational distribution one can find in the optimization procedure on the sampling data $\mathbb{D}$, assuming the parameters without prior assumptions are already optimized or fixed. We emphasize, that $Q_{\text{\tiny{PAC}}}$ is, depending on the choice of the model, of a special structure. If we talk about $Q_{\text{\tiny{PAC}}}$ in a specific context, $Q_{\text{\tiny{PAC}}}\in\mathcal{G}_{\hat{\Theta}}$ is the variational distribution of the specific model. Special instantiations of these can be found in beginning of Appendix~\ref{pacbayesian}.\\
Next, we introduce terms of loss functions in the context of Bayesian regression. Let $\mathbfcal{I}\in\mathbb{R}^\mathsf{K}$ be a Gaussian random vector with $\mathbfcal{I}\sim\mathcal{N}(\mathbf{0},\Sigma)$ and $E\in\mathbb{R}^{\mathsf{K}\times \mathsf{K}}$ symmetric, $\mathbf{e}\in\mathbb R^\mathsf{K}$, $e\in\mathbb R$ and a quadratic form defined as $\mathcal{Q}(\mathbfcal{I}) \stackrel{\mathrm{def}}= \mathbfcal{I}^TE\mathbfcal{I}+\mathbf{e}\mathbfcal{I}+e$. We say a loss $\ell$ is \textsl{sub-quadratic-form-Gaussian} if it can be described by a $\mathcal{Q}(\mathbfcal{I})$ random variable, i.e. its moment generating function is upper bounded by the one of a quadratic-form-Gaussian random variable $\mathcal{Q}(\mathbfcal{I})$, see Section 3.2 in~\cite{olkin1992quadratic}, which is
\begin{align*}
\vartheta(\lambda)& =\log\left(\underset{\substack{{\hat\theta\sim P_{\scalebox{.8}{\text{\tiny{PAC}}}}}, \boldsymbol{y}\sim  P_{\mathsf{K}}}}{\mathbf{E}}\left[e^{\lambda(\mathcal{L}_{P_\mathsf{K}}^{\ell}(\mathfrak{f}_{\uptheta})-\ell(\mathfrak{f}_{\uptheta},\mathbf{y}))}\right]\right)\\
&\stackrel{\mathrm{def}}{\leq} \log\left(\mathbf{E}\left[e^{\lambda \mathcal{Q}(\mathbfcal{I})}\right]\right)=-\frac{1}{2}\log\left(|I_{\mathsf{K}}-2\lambda E\Sigma|\right)+\frac{1}{2}\left(\lambda\mathbf{e}\right)^T(I_{\mathsf{K}}-2\lambda E\Sigma)^{-1}\Sigma\left(\lambda\mathbf{e}\right).
\end{align*}
In~\cite{germain2016pac} certain loss properties were relevant, which we will introduce as well. A loss $\ell$ is a \textsl{sub-Gaussian} loss, with a variance factor $s^2\in\mathbb R$, if it can be described by a Gaussian random variable $\mathcal{I}$, i.e. its moment generating function is upper bounded by the one of a Gaussian random variable of variance $s^2$. Furthermore, a loss $\ell$ is \textsl{sub-Gamma}, with a variance factor $s^2$ and scale $c$, if it can be described by a Gamma random variable $\mathcal{I}$, see also Section 2.3 in~\cite{boucheron2013concentration}.\\
For linear regression with simple priors, as well as simple data-distributions, these properties are enough. In our context of a DGP, more complex priors and a multi-dimensional data-distribution, it is natural to extend the above terms to the property sub-quadratic-form-Gaussian. Overall, when restricting to loss functions satisfying this new loss property in the context of DGPs, we expect tighter bound values then using the sub-Gaussian or sub-Gamma loss property. The reason for this is, when calculating the explicit expression for the moment generating function, further estimations would be necessary.

\section{The Deep Recurrent Gaussian Process models}
\label{sec:The Deep Gaussian Process models}
In this section, we shortly introduce two of the four DGP models which we already mentioned in Section~\ref{sec:introduction} and end of Section~\ref{sec:RelatedWork}, in particular for modeling time series data. These two models are compared in our experiments for the specific PAC-Bayesian statements. We follow \cite{mattos2015recurrent, foell2019deep} and refer there for details. The detailed structure of the DRGP with $L+1$ GP layers, where $\mathpzc{f}^{(l)}$ is a GP , is given by\vspace{-0.25cm}
\begin{align*}
\quad&\mathbf{h}^{i,(l)} = \mathrm{f}^{(l)}(\mathbf{X}^{(l)})+\boldsymbol\epsilon_i^
{\boldsymbol{h}^{(l)}},&&\text{with prior}\quad \boldsymbol{f}^{(l)}_{\mathbf{X}^{(l)}}
\sim\mathcal{N}(\mathbf{0},K_{\mathsf{K}\mathsf{K}}^{(l)}), && l = 1,\numberthis\label{gps1}
\dots,L\\
\quad&\mathbf{y}^i = \mathrm{f}^{(l)}(\mathbf{X}^{(l)})+\boldsymbol\epsilon_i^
{\boldsymbol{y}},&&\text{with prior}\quad \boldsymbol{f}^{(l)}_{\mathbf{X}^{(l)}}
\sim\mathcal{N}(\mathbf{0},K_{\mathsf{K}\mathsf{K}}^{(l)}), && l = L+1,\numberthis\label{gps2}
\end{align*}
\vspace{-0.25cm}
\\
with $\boldsymbol\epsilon_i^{\boldsymbol{h}^{(l)}}\sim \mathcal{N}(0,(\upsigma_{\text{noise}}
^{(l)})^2I_{\mathsf{K}})$, for $i = 1,\dots,\bar{N}$. We also write $\mathbf{h}^{i,(l)} \stackrel{\mathrm{def}}= [\mathrm{h}_1^{i,(l)},\dots,\mathrm{h}_{\mathsf{K}}^{i,(l)}]^T$ for $l=1,\dots,L+1$, $\mathbf{h}^{i,(L+1)} \stackrel{\mathrm{def}}= \mathbf{y}^{i}$ and respective the random vectors $\boldsymbol{h}^{(l)}$, $\boldsymbol{y}$. Here the measurements from $\mathsf{k}=1,\dots,\mathsf{K}$ are of time order. The matrix $K_{\mathsf{K}\mathsf{K}}$ represents a covariance matrix for a given covariance function $k^{(l)}$ and a set  of input-data $\mathbf{X}^{(l)}=[\mathbf{x}_{1}^{(l)},\dots,\mathbf{x}_{\mathsf{K}}^{(l)}]^T$, again of time order, for fixed time horizons $H_{\mathbf{x}}$, $H_{\mathrm{h}}$, is specified as\vspace{0.2cm}
\begin{flalign*}[left={\mathbf{x}_\mathsf{k}^{(l)} \stackrel{\mathrm{def}}
= \empheqlbrace}]
     && \begin{bmatrix}\mathbf{h}_{\mathsf{k}-1}^{(1)},\mathbf{\bar{x}}
_{\mathsf{k}-1}\end{bmatrix}^T&&\kern-1em\stackrel{\mathrm{def}}=&\left[\left[\mathrm{h
}_{\mathsf{k}-1}^{(1)},\dots,\mathrm{h}_{\mathsf{k}-H_{\mathrm{h}}}^{(1)}\right],\left[\mathbf{x
}_{\mathsf{k}-1},\dots,\mathbf{x}_{\mathsf{k}-H_{\mathbf{x}}}\right]\right]^T,&& l = 1
\numberthis\label{INPUT}\\
     && \quad\begin{bmatrix}\mathbf{h}_{\mathsf{k}-1}^{(l)},
\mathbf{h}_{\mathsf{k}}^{(l-1)}\end{bmatrix}^T&&\kern-1em\stackrel{\mathrm{def}}=&
\left[\left[\mathrm{h}_{\mathsf{k}-1}^{(l)},\dots,\mathrm{h}_{\mathsf{k}-H_{\mathrm{h}}}^{(l)}
\right],\left[\mathrm{h}_{\mathsf{k}}^{(l-1)},\dots,\mathrm{h}_{\mathsf{k}-H_{\mathrm{h}} + 1}^{
(l-1)}\right]\right]^T,&&l = 2,\dots,L\\
     && \mathbf{h}_{\mathsf{k}}^{(L)}&&\kern-1em\stackrel{\mathrm{def}}=&
\left[\mathrm{h}_{\mathsf{k}}^{(L)},\dots,\mathrm{h}_{\mathsf{k}-H_{\mathrm{h}} + 1}^{(L)}\right
]^T,&&l = L + 1,
\end{flalign*}
where $\mathbf{x}_\mathsf{k}^{(1)}\in\mathbb R^{H_{\mathrm{h}}+H_{\mathbf{x}}Q_{\mathbf{x}}}$, $\mathbf{x}_\mathsf{k}^{(l)}\in\mathbb R^{2H_h}$ for $l = 2,\dots,L$, $\mathbf{x}_\mathsf{k}^{(L+1)}\in\mathbb R^{H_{\mathrm{h}}}$, for $\mathsf{k} = 1,\dots,\mathsf{K}$. The $\mathsf{K}$ different state variables in this model are represented by $\mathbf{x}_\mathsf{k}^{(1)}$, more precisely by $\bar{\mathbf{x}}_\mathsf{k}$, for $\mathsf{k} = 1,\dots,\mathsf{K}$. For the other variables $\mathbf{x}_\mathsf{k}^{(1)}$, we refer to them as pseudo-states. The iteration until $\bar{N}$ just represents the amount of different output observation (one observation of a $\mathsf{K}$-dimensional vector) available for training. Later on, when variationally approximating the $\mathrm{h}_{\mathsf{k}}^{i,(l)}$, we derive for all input-data $\mathbf{X}^{(l)}$ the amount of $\mathsf{K}$ variational states, as we have $\mathsf{K}$ variational mean and variance parameters for these latent output-data points. That means, we are intended to model the output-data $\mathbf{Y}\in\mathbb R^{\mathsf{K}\times \bar{N}}$ on $\mathsf{K}$ variational states for all GP layers. This comes natural, when the unknown multivariate data-distribution $P_{\mathsf{K}}$ is of size $\mathsf{K}$. Be aware of that we are not restricted to these states neither in practice nor in theory (the posterior predictive stochastic process is well defined over the whole axis). Because of the independence assumption, the actual $K_{\hat{N}\hat{N}}^{(l)}$ covariance matrix, where $\hat{N}=\mathsf{K}\bar{N}$, collapses to a size of $\mathsf{K}\times \mathsf{K}$ in the modeling task.\\
Depending on the sparsity assumptions, we have different sparse covariance functions, as well as different priors $P_{\text{\scalebox{.8}{\tiny{PAC}}}}=P_{\text{\scalebox{.8}{\tiny{REV}}}}$ and variational distributions $Q_{\text{\scalebox{.8}{\tiny{PAC}}}}=Q_{\text{\scalebox{.8}{\tiny{REV}}}}$. In this paper, the sparse variational framework of~\cite{mattos2015recurrent} is called REVARB-Nyström and of~\cite{foell2019deep} REVARB-(V)SS. We refer to Appendix~\ref{pacbayesian} and~\cite[Section 4.]{mattos2015recurrent} for details. Our theorems generalize over the specific sparse variational framework, therefore we do not specify a specific sparse covariance function from the beginning. 

\section{PAC-Bayesian Bounds for DGPs}
\label{sec:PACBayesianBoundsforDGPs}
Our theorems hold for the DGPs of~\cite{damianou2015deep, mattos2015recurrent, cutajar2016practical, foell2019deep}. We derive these explicitly for the case of~\cite{foell2019deep} in the Appendix and also show how they can be adapted to~\cite{mattos2015recurrent}. Be aware of, that the inequality statements which follow have two-sided versions, which means they hold for the absolute value $|\cdot|$, see Appendix~\ref{definition}.
We state here the Theorem 4.1 from~\cite{alquier2015properties} for the case of empirical bounds, which we refine for our case.\\
\\
\noindent{\bf Theorem 1. (Empirical Bound)} [\cite{alquier2015properties}] {\it Given a data distribution $ P_\mathsf{K}$, a hypothesis set $\mathcal F_{\Theta}$, a loss function $\ell:\mathcal F_{\Theta}\times\mathbfcal{Y}\to\mathbb R$, a set of distributions $\mathcal{G}_{\hat{\Theta}}$, a prior $P_{\text{\scalebox{.8}{\tiny{PAC}}}}$ in $\mathcal{G}_{\hat{\Theta}}$, a posterior $Q_{\text{\scalebox{.8}{\tiny{PAC}}}}$ in $\mathcal{G}_{\hat{\Theta}}$, a $\delta\in (0,1]$, then with probability at least $1-\uptau$ over $\mathbbmsl{D}\sim (P_\mathsf{K})^N$ we have for all $Q_{\text{\scalebox{.8}{\tiny{PAC}}}}\text{ in }\mathcal{G}_{\hat{\Theta}}$:
\begin{align*} 
&\underset{{\hat\theta\sim Q_{\text{\scalebox{.8}{\tiny{PAC}}}}}}{\mathbf{E}}\left[\mathcal{L}_{P_\mathsf{K}}^{\ell}(\mathfrak{f}_{\uptheta})\right]-\underset{{\hat\theta\sim Q_{\text{\scalebox{.8}{\tiny{PAC}}}}}}{\mathbf{E}}\left[\mathcal{L}_{\mathbb{D}}^{\ell}(\mathfrak{f}_{\uptheta})\right]\leq \frac{1}{\lambda}\left(\mathbf{KL}(Q_{\text{\scalebox{.8}{\tiny{PAC}}}}||P_{\text{\scalebox{.8}{\tiny{PAC}}}})+\log\left(\frac{1}{\uptau}\right)+\Psi^{\ell}(\lambda,N)\right),
\end{align*}
\vspace{-0.3cm}
\begin{align*}
&\text{ where }&\Psi^{\ell}(\lambda,N)=\log\left(\underset{\substack{{\hat\theta \sim P_{\scalebox{.8}{\text{\tiny{PAC}}}}},{\mathbbmsl{D}'\sim (P_\mathsf{K})^N}}}{\mathbf{E}}\left[e^{\lambda(\mathcal{L}_{P_\mathsf{K}}^{\ell}(\mathfrak{f}_{\scalebox{.6}{$\uptheta$}})-\mathcal{L}_{\mathbb{D}'}^{\ell}(\mathfrak{f}_{\scalebox{.6}{$\uptheta$}}))}\right]\right).
\end{align*}
}The main difficulty to derive explicit expressions for the right hand side of the inequality for DGPs is calculating the expression $\Psi^{\ell}(\lambda,N)$ explicitly. Given the likelihood of the DGPs, as well as the priors and variational distributions, this is not straightforward but possible under some specific assumptions, as the likelihood decomposes into a sum of $L+1$ terms representing the likelihood of simple sparse GPs, see Appendix \ref{pacbayesian}, Equation \ref{empirical}. These assumptions are:\\
\\
\textbf{Assumption 1 (A1) [\text{Stochastic Lipschitz condition}]:}
\begin{flalign*}
& \text{For functions $\mathfrak{v}^{(l)}:\mathbb R^{p^{(l)}+v^{(l)}}\to\mathbb R$, $\mathfrak{f}:\mathbb R^{p^{(l)}}\to\mathbb R$, where one is a stochastic version of the other}\\
& \text{(some parameters in the function $\mathfrak{f}^{(l)}$ are assumed to have a mean $\uptheta_\mathrm{m}^{(l)}$, variance $\uptheta_\mathrm{v}^{(l)}$, which results}\\
& \text{in the function $\mathfrak{v}^{(l)}$), we assume the property that $\exists \mathbf{S}\in\mathbb R_\geq$ such that $\forall (\uptheta^{(l)},\mathbf{0})^T,(\check\uptheta^{(l)},\uptheta_\mathrm{m}^{(l)},\uptheta_\mathrm{v}^{(l)})^T$}:\\
&\quad\quad\frac{\Vert\mathfrak{f}^{(l)}(\uptheta^{(l)})-\mathfrak{v}^{(l)}(\check\uptheta^{(l)},\uptheta_\mathrm{m}^{(l)},\uptheta_\mathrm{v}^{(l)})\Vert^2}{\Vert(\uptheta^{(l)},\mathbf{0})^T-(\check\uptheta^{(l)},\uptheta_\mathrm{m}^{(l)},\uptheta_\mathrm{v}^{(l)})^T\Vert^2} \leq \mathbf{S}^2,\quad \Vert(\uptheta^{(l)},\mathbf{0})^T-(\check\uptheta^{(l)},\uptheta_\mathrm{m}^{(l)},\uptheta_\mathrm{v}^{(l)})^T\Vert^2\neq 0.
\end{flalign*}
\\
\textbf{Assumption 2 (A2) [\text{Fubinis Theorem}]:}
\begin{flalign*}
& \underset{\substack{{\hat\theta \sim P_{\scalebox{.8}{\text{\tiny{REV}}}}},{ \boldsymbol{y}\sim P_\mathsf{K}}}}{\mathbf{E}}[-\log(p(\mathbf{y}|\uptheta,\mathbf{X}))] = \underset{\substack{{\boldsymbol{y}\sim P_\mathsf{K}},{ \hat\theta \sim P_{\scalebox{.8}{\text{\tiny{REV}}}}}}}{\mathbf{E}}[-\log(p(\mathbf{y}|\uptheta,\mathbf{X}))],\\
 &\underset{\substack{{\hat\theta \sim Q_{\scalebox{.8}{\text{\tiny{REV}}}}},{ \boldsymbol{y}\sim P_\mathsf{K}}}}{\mathbf{E}}[-\log(p(\mathbf{y}|\uptheta,\mathbf{X}))] = \underset{\substack{{\boldsymbol{y}\sim P_\mathsf{K}},{ \hat\theta \sim Q_{\scalebox{.8}{\text{\tiny{REV}}}}}}}{\mathbf{E}}[-\log(p(\mathbf{y}|\uptheta,\mathbf{X}))].
\end{flalign*}
\\
\textbf{Assumption 3 (A3) [\text{point-wise bounded Covariance matrix}]:}
\begin{flalign*}
&\textbf{Cov}\left[\boldsymbol{h}^{(l)}\right]\text{ for $l=1,\dots,L+1$ GP layers is bounded for each component .}
\end{flalign*}
Assumption 1 evolves naturally, when we are dealing with variational approximations of predefined models. It just expresses, that two functions with the same output-space and similar inputs should hold the property of the simple Lipschitz condition. Here, we identify the mean and variationally approximated value with each other, and stack the input space for the function $\mathfrak{f}$ with variance $\boldsymbol{0}$. For our case of DGPs, we use this assumption for each GP layers' predictor $\mathfrak{f}^{(l)}$, for $l=1\dots,L+1$. For the first GP layer, the two functions are indexed with the fixed input-states $\bar{\mathbf{x}}_\mathsf{k}$. Together with the bound on the states input-space, which we assumed at the beginning, this extends the stochastic Lipschitz condition for the first layer by considering the states as new inputs. When we bound the functions input-space, this condition is equivalent to bounding the functions output-space (again, be aware of that the predictor $\mathfrak{f}^{(l)}$ is the mean-function, not the samples of the GPs - the samples are \textit{not} bounded). Following~\cite{NIPS20177100}, Assumption 2 is a relative mild assumption and we refer to their paper for further cases, when this holds. Assumption 3 is obvious and no further explanation is needed.\\
In~\cite{alquier2015properties, germain2016pac} properties of variational approximations and a connection between PAC-Bayesian theory and Bayesian Inference in terms of the marginal likelihood were shown, which we can adapt for our case. We see, that we can have the same intuition about this link of PAC-Bayesian bound and the marginal likelihood for our models, which use a variational bound instead of the marginal likelihood. We show this link between our variational bounds $\mathbfcal{L}_{\text{\scalebox{.8}{\tiny{REV}}}}$ and the PAC-Bayesian generalization risk bound choosing $\lambda=N$. We derive with the negative likelihood loss function $\ell=\ell_{\text{nll}}$ and $P_{\text{\scalebox{.8}{\tiny{PAC}}}}=P_{\text{\scalebox{.8}{\tiny{REV}}}}$, $Q_{\text{\scalebox{.8}{\tiny{PAC}}}}=Q_{\text{\scalebox{.8}{\tiny{REV}}}}$ the priors and variational distributions of the REVARB-frameworks:\vspace{0.2cm}
\begin{flalign*} 
N\underset{{\hat\theta\sim Q_{\text{\tiny{PAC}}}}}{\mathbf{E}}\left[\mathcal{L}_{\mathbb{D}}^{\ell_{\text{nll}}}(\mathfrak{f}_{\uptheta})\right] + \mathbf{KL}(Q_{\text{\scalebox{.8}{\tiny{PAC}}}}||P_{\text{\scalebox{.8}{\tiny{PAC}}}})& =N\underset{{\hat\theta\sim Q_{\text{\scalebox{.8}{\tiny{REV}}}}}}{\mathbf{E}}\left[\mathcal{L}_{\mathbb{D}}^{\ell_{\text{nll}}}(\mathfrak{f}_{\uptheta})\right]+ \mathbf{KL}(Q_{\text{\scalebox{.8}{\tiny{REV}}}}||P_{\text{\scalebox{.8}{\tiny{REV}}}})\numberthis\label{PACVAR}\\
&=\underset{{\hat\theta\sim Q_{\text{\scalebox{.8}{\tiny{REV}}}}}}{\mathbf{E}}\left[-\sum\nolimits_{i=1}^N\log(p(\mathbf{y}^i|\uptheta,\mathbf{X}))\right]+ \mathbf{KL}(Q_{\text{\scalebox{.8}{\tiny{REV}}}}||P_{\text{\scalebox{.8}{\tiny{REV}}}})\\
&=\underset{{\hat\theta\sim Q_{\text{\scalebox{.8}{\tiny{REV}}}}}}{\mathbf{E}}\left[-\log\left(\prod\nolimits_{i=1}^N p(\mathbf{y}^i|\uptheta,\mathbf{X})\right)\right]+ \mathbf{KL}(Q_{\text{\scalebox{.8}{\tiny{REV}}}}||P_{\text{\scalebox{.8}{\tiny{REV}}}})\\
&=-\underset{{\hat\theta\sim Q_{\text{\scalebox{.8}{\tiny{REV}}}}}}{\mathbf{E}}[\mathcal{G}_{\text{\scalebox{.8}{\tiny{REV}}}}] + \mathbf{KL}(Q_{\text{\scalebox{.8}{\tiny{REV}}}}||P_{\text{\scalebox{.8}{\tiny{REV}}}})=-\mathbfcal{L}_{\text{\scalebox{.8}{\tiny{REV}}}}.
\end{flalign*}
$\mathbfcal{L}_{\text{\scalebox{.8}{\tiny{REV}}}}$ can e.g. be found in~\cite{foell2019deep}, Section 4.2, Equation (14), (15).
We now can state the following theorems, which themselves directly follow from Theorem 1, the above Equation \eqref{PACVAR} and our introduced new term for loss function (Proof in Appendix~\ref{pacbayesian}). It is similar to the Corollary 5 of~\cite{germain2016pac} but with the extended term for the loss function, which is more appropriate and which shows not just a connection to the variational bound for the DGP models, but even convergence in the sense of consistency with minimal assumptions.\\
\\\newpage
\noindent{\bf Theorem 2. (Empirical Bound 1 for DGPs \eqref{gps1}, \eqref{gps2}, \eqref{INPUT}, \eqref{PACVAR})} {\it Given a data distribution $P_\mathsf{K}$, a hypothesis set $\mathcal F_{\Theta}$, the loss $\ell=\ell_{\text{nll}}$ which is quadratic-form-Gaussian, associated priors $P_{\text{\scalebox{.8}{\tiny{PAC}}}}=P_{\text{\scalebox{.8}{\tiny{REV}}}}$ in $\mathcal{G}_{\hat{\Theta}}$, a posterior $Q_{\text{\scalebox{.8}{\tiny{PAC}}}}=Q_{\text{\scalebox{.8}{\tiny{REV}}}}$  in $\mathcal{G}_{\hat{\Theta}}$, a $\uptau\in (0,1]$, Assumption 1: a stochastic Lipschitz condition for the raw model and the variational model, Assumption 2: Fubinis theorems, Assumption 3: point-wise bounded $\textbf{Cov}\left[\boldsymbol{h}^{(l)}\right]$, bounded input-space for the mean-function for all $\mathsf{k}=1,\dots,\mathsf{K}$, $l=1,\dots,L+1$, then with probability at least $1-\uptau$ over $\mathbbmsl{D}\sim ( P_\mathsf{K})^N$ we have for all $Q_{\text{\scalebox{.8}{\tiny{REV}}}}\text{ in }\mathcal{G}_{\hat{\Theta}}$ with $\lambda=N$:
\begin{align*} 
\underset{{\hat\theta\sim Q_{\text{\scalebox{.8}{\tiny{REV}}}}}}{\mathbf{E}}\left[\mathcal{L}_{ P_\mathsf{K}}^{\ell_{\text{nll}}}(\mathfrak{f}_{\uptheta})\right]\leq\frac{\mathcal{L}(N)}{N}-\frac{\log(e^{\mathbfcal{L}_{\text{\scalebox{.8}{\tiny{REV}}}}}\uptau)}{N},
\end{align*}
where we can choose $Q_{\text{\scalebox{.8}{\tiny{REV}}}}=Q_{\lambda}$ and $\mathcal{L}(N)$ is defined in the Appendix \ref{pacbayesian}, Equation \eqref{L2}.}\\
\\
\noindent{\bf Theorem 3. (Empirical Bound 2, Consistency for DGPs \eqref{gps1}, \eqref{gps2}, \eqref{INPUT})} {\it Given the same assumptions as in Theorem 2., then we have with $\lambda=\sqrt{N}$:
\begin{align*} 
\underset{{\hat\theta\sim Q_{\text{\scalebox{.8}{\tiny{REV}}}}}}{\mathbf{E}}\left[\mathcal{L}_{ P_\mathsf{K}}^{\ell_{\text{nll}}}(\mathfrak{f}_{\uptheta})\right]-\underset{{\hat\theta\sim Q_{\text{\scalebox{.8}{\tiny{REV}}}}}}{\mathbf{E}}\left[\mathcal{L}_{\mathbb{D}}^{\ell_{\text{nll}}}(\mathfrak{f}_{\uptheta})\right]\leq
\frac{\mathbf{KL}(Q_{\text{\scalebox{.8}{\tiny{REV}}}}||P_{\text{\scalebox{.8}{\tiny{REV}}}})+\log\left(\frac{1}{\uptau}\right)+\mathcal{L}(\sqrt{N})}{\sqrt{N}},
\end{align*}
where we can choose $Q_{\text{\scalebox{.8}{\tiny{REV}}}}=Q_{\sqrt{N}}$ and $\mathcal{L}(\sqrt{N})$ is defined in the Appendix \ref{pacbayesian}, Equation \eqref{L2}. 
We have convergence of order $\mathcal{O}\left(\frac{1}{\sqrt{N}}\right)$ to zero (consistency).
}\\
\\
Let $\mathbf{S}^{(l)}$ be the stochastic Lipschitz constants, $\updelta^{(l)}$ the bound on the input spaces, ${\displaystyle 1\!\!1_\mathsf{K}^T}\textbf{Var}\left[\boldsymbol{h}^{(l)}\right]$ the sum of the variance of the posterior predictive distribution of our variational models, where we have for \cite{foell2019deep} DRGP-SS explicitly:
\begin{align*}
{\displaystyle 1\!\!1_\mathsf{K}^T}\textbf{Var}\left[\boldsymbol{h}^{(l)}\right] &= (\mathbf{m}^{(l)}_*)^T(\Psi^{(l)}_{2,*}-(\Psi^{(l)}_{1,*})^T\Psi^{(l)}_{1,*})\mathbf{m}^{(l)}_*+\text{tr}\left(\Psi^{(l)}_{2,*}\mathbf{s}^{(l)}_*\right)+\mathsf{K}(\sigma_{\text{noise}*}^{(l)})^2,\\
\textbf{Cov}\left[\boldsymbol{h}^{(l)}\right] &= \Psi^{(l)}_{1,*}\mathbf{s}^{(l)}_*(\Psi^{(l)}_{1,*})^T\in\mathbb R^{\mathsf{K}\times \mathsf{K}}, \mathsf{k}\neq \hat{\mathsf{k}}, \mathsf{k},\hat{\mathsf{k}}=1,\dots,\mathsf{K}.
\end{align*}
For the variables involved in ${\displaystyle 1\!\!1_\mathsf{K}^T}\textbf{Var}\left[\boldsymbol{h}^{(l)}\right]$ and $\textbf{Cov}\left[\boldsymbol{h}^{(l)}\right]$ we refer to beginning of Appendix \ref{pacbayesian}.\\
The expression $\frac{\mathcal{L}(\sqrt{N})}{\sqrt{N}}$ for the DGPs then has the form:
{
\begin{align*}
&\sum\nolimits_{l=1}^{L+1}\frac{{\displaystyle 1\!\!1_\mathsf{K}^T}\textbf{Var}\left[\boldsymbol{h}^{(l)}\right]}{2(\sigma_{\text{noise}}^{(l)})^2}-\frac{\sqrt{N}}{2}\log\left(|I_\mathsf{K}+\frac{
 \textbf{Cov}\left[\boldsymbol{h}^{(l)}\right]}{\sqrt{N}(\sigma_{\text{noise}}^{(l)})^2}|\right)+\frac{1}{2\sqrt{N}}\left(\left[\frac{\mathbf{S}^{(l)}}{\updelta}\frac{1}{(\sigma_{\text{noise}}^{(l)})^2}\right]_{\mathsf{k}=1}^\mathsf{K}\right)^T\\
&\left(\textbf{Cov}\left[\boldsymbol{h}^{(l)}\right]\right)\left(I_\mathsf{K}+\frac{\textbf{Cov}\left[\boldsymbol{h}^{(l)}\right]}{\sqrt{N}(\sigma_{\text{noise}}^{(l)})^2}\right)^{-1}\left[\frac{\mathbf{S}^{(l)}}{\updelta}\frac{1}{(\sigma_{\text{noise}}^{(l)})^2}\right]_{\mathsf{k}=1}^\mathsf{K},
\end{align*}}%
We state next the Theorem 4.2 from~\cite{alquier2015properties} for the case of oracle-type inequalities, which we also refine for our case.\\
\\\newpage
\noindent{\bf Theorem 4.} (Oracle-type Bound) {\it Given a data distribution $ P_\mathsf{K}$, a hypothesis set $\mathcal F_{\Theta}$, a loss function $\ell:\mathcal F_{\boldsymbol\Theta}\times\mathbfcal{Y}\to\mathbb R$, a set of distributions $\mathcal{G}_{\hat{\Theta}}$, a prior $P_{\text{\scalebox{.8}{\tiny{PAC}}}}$ in $\mathcal{G}_{\hat{\Theta}}$, the best variational posterior $Q_{\lambda}$ in $\mathcal{G}_{\hat{\Theta}}$, a $\uptau\in (0,1]$, a real number $\lambda>0$, then with probability at least $1-\uptau$ over $\mathbbmsl{D}\sim ( P_\mathsf{K})^N$, we have for all $Q_{\text{\scalebox{.8}{\tiny{PAC}}}}$ in $\mathcal{G}_{\hat{\Theta}}$:
\begin{align*} 
\underset{{\hat\theta\sim Q_{\lambda}}}{\mathbf{E}}\left[\mathcal{L}_{P_\mathsf{K}}^{\ell}(\mathfrak{f}_{\uptheta})\right]\leq\inf\limits_{Q_{\text{\scalebox{.8}{\tiny{PAC}}}}\text{ in }\mathcal{G}_{\hat{\Theta}}}\left(\underset{{\hat\theta\sim Q_{\text{\scalebox{.8}{\tiny{PAC}}}}}}{\mathbf{E}}\left[\mathcal{L}_{P_\mathsf{K}}^{\ell}(\mathfrak{f}_{\uptheta})\right]+\frac{\mathbf{KL}(Q_{\text{\scalebox{.8}{\tiny{PAC}}}}||P_{\text{\scalebox{.8}{\tiny{PAC}}}})+\log\left(\frac{1}{\uptau}\right)+\Psi^{\ell}(\lambda,N)}{\lambda}\right).
\end{align*}}
\\
Theorem 4 is the point of view developed in~\cite{catoni2004statistical, catoni2007pac} and~\cite{dalalyan2008aggregation}. This statement enables us to compare our optimized variational error against the best possible aggregation procedure in $\mathcal {G}_{\hat{\Theta}}$ and the raw model error (no variational expectation).  
We now can state the following theorem, which itself directly follows from Theorem 4 and our introduced new property of loss function (Proof in Appendix \ref{pacbayesian2}). It shows an interesting bound, where the variance of the variational posterior predictive distribution and the variational variance of our variational approximation for the pseudo output-data $\mathbf{h}^{(l)}$ is involved. This statement gives insights in the convergence between the raw model error and our variational models. For details regarding the terms involved we again refer to the Appendix \ref{pacbayesian2} and \cite{foell2019deep}.\\
\\
\noindent{\bf Theorem 5.} (Oracle-type inequality) {\it Given a data distribution $P_\mathsf{K}$, a hypothesis set $\mathcal F_{\Theta}$, the loss $\ell=\ell_{\text{nll}}$ which is quadratic-form-Gaussian, associated priors $P_{\text{\scalebox{.8}{\tiny{PAC}}}}=P_{\text{\scalebox{.8}{\tiny{REV}}}}$ in $\mathcal{G}_{\hat{\Theta}}$, the best variational posterior $Q_{\lambda}$ in $\mathcal{G}_{\hat{\Theta}}$, a $\uptau\in (0,1]$, Assumption 1: a stochastic Lipschitz condition for the raw model and the variational model, Assumption 2: Fubinis theorems, Assumption 3: point-wise bounded $\textbf{Cov}\left[\boldsymbol{h}^{(l)}\right]$, bounded input-space for the mean-function for all $\mathsf{k}=1,\dots,\mathsf{K}$, $l=1,\dots,L+1$, then with probability at least $1-\uptau$ over $\mathbbmsl{D}\sim ( P_\mathsf{K})^N$ we have for all $Q_{\text{\scalebox{.8}{\tiny{REV}}}}\text{ in }\mathcal{G}_{\hat{\Theta}}$ with $\lambda=\sqrt{N}$:
\begin{align*}
&\underset{{\hat\theta\sim Q_{\sqrt{N}}}}{\mathbf{E}}\left[\mathcal{L}_{P_\mathsf{K}}^{\ell_{\text{nll}}}(\mathfrak{f}_{\uptheta})\right]- \mathcal{L}_{P_\mathsf{K}}^{\ell_{\text{nll}}}(\mathfrak{f}_{\uptheta_*})\leq \mathcal{L}_{\text{\scalebox{.8}{\tiny{Ora}}}}+\frac{\mathbf{KL}(Q_{*}||P_{\text{\scalebox{.8}{\tiny{REV}}}})+\log\left(\frac{1}{\uptau}\right)+\mathcal{L}(\sqrt{N})}{\sqrt{N}},
\end{align*}
where $\mathcal{L}_{\text{\scalebox{.8}{\tiny{Ora}}}}=\sum\nolimits_{l=1}^{L+1}\frac{{\displaystyle 1\!\!1_\mathsf{K}^T}\left(\boldsymbol\uplambda_*^{(l)}+\textbf{Var}\left[\boldsymbol{f}^{(l)}_{\mathbf{X}^{(l)},\uptheta_*}\right]\right)}{2(\sigma_{\text{noise}*}^{(l)})^2}$.
We have convergence of order $\mathcal{O}\left(\frac{1}{\sqrt{N}}\right)$ to $\mathcal{L}_{\text{\scalebox{.8}{\tiny{Ora}}}}$.
}\\
\\
Next, we make use of the union bound (UB) and some upcoming properties to additionally derive for all theorems new forms. The goal is to derive statements, which explicitly involve the first GP layers' input-dimension. For the task of modeling time-series data, the input dimension grows linearly with the time-horizons, hence $Q$ is of big interest. Therefore, we introduce the term of the smallest covering numbers of a metric space, here the input-space $\mathbb R^Q$ with the Euclidean metric. Let $A\subset \mathbb R^Q$, then for $\epsilon>0$ and $B_{\|.\|}(\mathbf{x},\epsilon)\stackrel{\mathrm{def}}{=}\{\mathbf{x}'\in\mathbb R^Q: \begin{Vmatrix}\mathbf{x}-\mathbf{x}'\end{Vmatrix}\leq \epsilon\}$, we define
\begin{align*}
\mathcal{C}(A,\|.\|,\epsilon)\stackrel{\mathrm{def}}{=}\text{min}\{M\geq 1:\exists\;\mathbf{a}_1,\dots,\mathbf{a}_M\in\mathbb R^Q, A\subset \bigcup\nolimits_{i=1}^M B_{\|.\|}(\mathbf{a}_i,\epsilon)\},
\end{align*}
Let $\mathfrak{f}_{\uptheta}(\mathbf{x}^{(1)})$ be the predictor in form of the composition of the GP predictors $\mathfrak{f}^{(1)},\dots,\mathfrak{f}^{(L+1)}$ in the Appendix \ref{pacbayesian}, Equations (\ref{predictor1})-(\ref{predictor3}) for a new input $\mathbf{x}^{(1)}$ following the notation in Section \ref{sec:The Deep Gaussian Process models} and $\mathcal{K}\stackrel{\mathrm{def}}=\{\mathfrak{f}_{\uptheta}(\mathbf{x}^{(1)}):\mathbf{x}^{(1)}\in rB_{\|.\|},r\in\mathbb R\}$, where $B_{\|.\|}=[-1,1]^Q$ with $Q$ the input-dimension.\\
Now, we let $\mathbf{x}_1^{(1)},\dots,\mathbf{x}_\mathsf{K}^{(1)}\in rB_{\|.\|}$ be chosen in such a way, that $\bigcup\nolimits_{\mathsf{k}=1}^{\mathsf{K}}B_{\|.\|}(\mathfrak{f}_{\uptheta}(\mathbf{x}_{\mathsf{k}}^{(1)}),\epsilon)$ covers $\mathcal{K}$ w.r.t. $\|.\|$ with $\mathsf{K}=\mathcal{C}(\mathcal{K},\|.\|,\epsilon)$, then we get with  $R=\frac{\mathbf{KL}(Q_{\text{\scalebox{.8}{\tiny{REV}}}}||P_{\text{\scalebox{.8}{\tiny{REV}}}})+\log\left(\frac{1}{\uptau}\right)+\mathcal{L}\left(\lambda\right)}{\lambda}$, the union bound, $\mathfrak{f}_{{\mathsf{k}},\uptheta}$ the predictor for the ${\mathsf{k}}$-th state, Appendix~\ref{extensions} (AE) and Theorem 3 as show-case:
\begin{flalign*} 
(P_{\mathsf{K}})^N&\left(\underset{{\theta\sim Q_{\lambda}}}{\mathbf{E}}\left[\mathcal{L}_{P_{\mathsf{K}}}^{\ell_{\text{nll}}}(\mathfrak{f}_{\uptheta})\right]-\underset{{\theta\sim Q_{\lambda}}}{\mathbf{E}}\left[\mathcal{L}_{\mathbb{D}}^{\ell_{\text{nll}}}(\mathfrak{f}_{\uptheta})\right]\leq R\right)\numberthis\label{union}\\
&\geq (P_{\mathsf{K}})^N\left(\mathsf{K}\sup\limits_{{\mathsf{k}}=1,\dots,\mathsf{K}}\left\{\underset{{\theta\sim Q_{\lambda}}}{\mathbf{E}}\left[\mathcal{L}_{P_{\mathsf{K}}}^{\ell_{\text{nll}}}(\mathfrak{f}_{{\mathsf{k}},\uptheta})\right]-\underset{{\theta\sim Q_{\lambda}}}{\mathbf{E}}\left[\mathcal{L}_{\mathbb{D}}^{\ell_{\text{nll}}}(\mathfrak{f}_{{\mathsf{k}},\uptheta})\right]\right\}\leq R\right)\\
&=(P_{\mathsf{K}})^N\left(\bigcap\limits_{{\mathsf{k}}=1,\dots,\mathsf{K}}\left\{\underset{{\theta\sim Q_{\lambda}}}{\mathbf{E}}\left[\mathcal{L}_{P_{\mathsf{K}}}^{\ell_{\text{nll}}}(\mathfrak{f}_{{\mathsf{k}},\uptheta})\right]-\underset{{\theta\sim Q_{\lambda}}}{\mathbf{E}}\left[\mathcal{L}_{\mathbb{D}}^{\ell_{\text{nll}}}(\mathfrak{f}_{{\mathsf{k}},\uptheta})\right]\leq \frac{R}{\mathsf{K}}\right\}\right)\\
&\stackrel{\mathrm{(UB)}}{\geq}\sum\nolimits_{\mathsf{k}=1}^\mathsf{K} (P_{\mathsf{K}})^N\left(\underset{{\theta\sim Q_{\lambda}}}{\mathbf{E}}\left[\mathcal{L}_{P_{\mathsf{K}}}^{\ell_{\text{nll}}}(\mathfrak{f}_{{\mathsf{k}},\uptheta})\right]-\underset{{\theta\sim Q_{\lambda}}}{\mathbf{E}}\left[\mathcal{L}_{\mathbb{D}}^{\ell_{\text{nll}}}(\mathfrak{f}_{{\mathsf{k}},\uptheta})\right]\leq R_{\lambda'=\lambda \mathsf{K}}\right)-(\mathsf{K}-1)\\
&\stackrel{\mathrm{(AE)}}{\geq} \sum\nolimits_{{\mathsf{k}}=1}^\mathsf{K} \left(1-\uptau\frac{\mathcal{L}_{\mathsf{k}}\left(\lambda'\right)}{\mathcal{L}\left(\lambda'\mathsf{K}^{-1}\right)}\right) - (\mathsf{K} -1)\\
&= 1-\mathsf{K}\uptau\frac{\mathcal{L}_{\mathsf{k}}\left(\lambda'\right)}{\mathcal{L}\left(\lambda'\mathsf{K}^{-1}\right)}.
\end{flalign*}
Furthermore, assuming that for every predictor $\mathfrak{f}^{(l)}$ of the GPs in the DGP the standard Lipschitz condition w.r.t. to the state input $\mathbf{x}_\mathsf{k}^{(l)}$ holds, where $\mathbf{L}^{(l)}$ is the Lipschitz constant, (this is automatically fullfilled for all layers with the stochastic Lipschitz condition) we obtain
\begin{align*}
\|\mathfrak{f}_{\uptheta}(\mathbf{x}_\mathsf{k}^{(1)})-\mathfrak{f}_{\uptheta}({\mathbf{x}}_\mathsf{k}'^{(1)})\|= \mathbf{L}^{(L+1)}\|\mathbf{x}_\mathsf{k}^{(L)}-{\mathbf{x}}_\mathsf{k}'^{(L)}\|=\dots=\left(\prod\nolimits_{l=1}^{L+1}\mathbf{L}^{(l)}\right)\|\mathbf{\bar{x}}_\mathsf{k}-\mathbf{\bar{x}_\mathsf{k}'}\|,
\end{align*}
where we assume that the auto-regressive values in Equation \eqref{INPUT} are equal and $H_{\mathrm{h}}=1$.\\
We define $\mathbf{L}=\prod\nolimits_{l=1}^{L+1}\mathbf{L}^{(l)}$, and with $\epsilon'=\frac{\epsilon}{r}$ we can now derive 
\begin{align*}
\mathcal{C}(\mathcal{K},\|.\|,\epsilon) &\leq\mathcal{C}(rB_{\|.\|},\|.\|,\frac{\epsilon}{\mathbf{L}})=\mathcal{C}(B_{\|.\|},\|.\|,\frac{\epsilon}{r\mathbf{L}})\leq \left(\frac{\epsilon'}{\mathbf{L}}\right)^{-Q}.\numberthis\label{cover}
\end{align*}
After a transformation $1-\uptau'=1-\mathsf{K}\uptau\frac{\mathcal{L}_{{\mathsf{k}}}\left(\lambda'\right)}{\mathcal{L}\left(\lambda'\mathsf{K}^{-1}\right)}$ for the inequality derived in~\eqref{union} we can now use this derived estimation \eqref{cover} and the definitions to derive new bounds by 
\begin{align*} 
&(P_{\mathsf{K}})^N\left(\underset{{\theta\sim Q_{\lambda}}}{\mathbf{E}}\left[\mathcal{L}_{P_{\mathsf{K}}}^{\ell_{\text{nll}}}(\mathfrak{f}_{\uptheta})\right]-\underset{{\theta\sim Q_{\lambda}}}{\mathbf{E}}\left[\mathcal{L}_{\mathbb{D}}^{\ell_{\text{nll}}}(\mathfrak{f}_{\uptheta})\right]\leq \hat{R}\right)\geq 1-\uptau',
\end{align*}
with
\begin{align*} 
\hat{R}=\frac{\mathbf{KL}(Q_{\text{\scalebox{.8}{\tiny{REV}}}}||P_{\text{\scalebox{.8}{\tiny{REV}}}})+\log\left(\frac{\mathsf{K}}{\uptau'}\right)+\log\Bigl(\frac{\mathcal{L}\left(\lambda'\mathsf{K}^{-1}\right)}{\mathcal{L}_{\mathsf{k}}\left(\lambda'\right)}\Bigr)+\mathcal{L}(\lambda)}{\lambda},
\end{align*}
and as $\log(\mathsf{K})\leq Q \log\left(\frac{\mathbf{L}}{\epsilon'}\right)$ and choosing $\epsilon' = \frac{1}{N}$ we come to 
\begin{align*} 
\hat{R}=\frac{\mathbf{KL}(Q_{\text{\scalebox{.8}{\tiny{REV}}}}||P_{\text{\scalebox{.8}{\tiny{REV}}}})+Q\log\left(N\mathbf{L}\right)+\log\left(\frac{1}{\uptau'}\right)+\log\Bigl(\frac{\mathcal{L}\left(\lambda'\mathsf{K}^{-1}\right)}{\mathcal{L}_{\mathsf{k}}\left(\lambda'\right)}\Bigr)+\mathcal{L}(\lambda)}{\lambda}.\numberthis\label{newbound}
\end{align*}
This new bound~(\ref{newbound}) holds for Theorem 3, but can be extended to the other ones in the same way. It is convenient to choose $\updelta = Q$ involved in Equation~(\ref{kappa}), $\mathcal{L}(\lambda)$ and the Appendix \ref{pacbayesian}, the input-space bound involved in the stochastic Lipschitz condition. Then we have a multiplicative impact of $Q$ in nearly all terms regarding $\lambda=N,\sqrt{N}\to\infty$. The expression $\log\Bigl(\frac{\mathcal{L}\left(\lambda'\mathsf{K}^{-1}\right)}{\mathcal{L}_{\mathsf{k}}\left(\lambda'\right)}\Bigr)$ converges of order $\mathcal{O}\left(\frac{1}{\sqrt{N}}\right)$ to a fixed value.

\subsection{Experimental results}
\label{experiments}
In our experiments we want to show the evolution of the error of the consistency result for $\lambda=\sqrt{N}\to\infty$, $\mathcal{L}(\lambda)$ of Theorem 3 for \cite{mattos2015recurrent} (RGP) and \cite{foell2019deep} (DRGP-(V)SS)  on the data-sets involved in~\cite{foell2019deep}. We choose $\uptau=0.5$, $\lambda=\sqrt{N}$ and we create, based on the model training result of the specific model from the training data-set $\mathbf{Y}\in\mathbb R^{\mathsf{K}\times \bar{N}}$, where $\bar{N}=1$, and the input-data is chosen as the $\mathsf{K}$ states, new measurements $\mathbf{y}^i\in\mathbb R^\mathsf{K}$, $i=1,\dots,N$, by adding noise $N=50000$ times to the predicted training output-data $\mathbf{y}^{\text{pred}}\in \mathbb R^\mathsf{K}$. The noise follows the model variance prediction on these states, which also comes from the training on the specific data-sets. In Figure \ref{fig:VIS1} in the Appendix~\ref{experimentsA} we see the predicted training output-data $\mathbf{y}^{\text{pred}}$ of the Actuator data-set and its created samples $\mathbf{y}^i$. 
\begin{figure}[H]
\vspace{-0.2cm}
\centering
	{\includegraphics[width=0.4\textwidth]{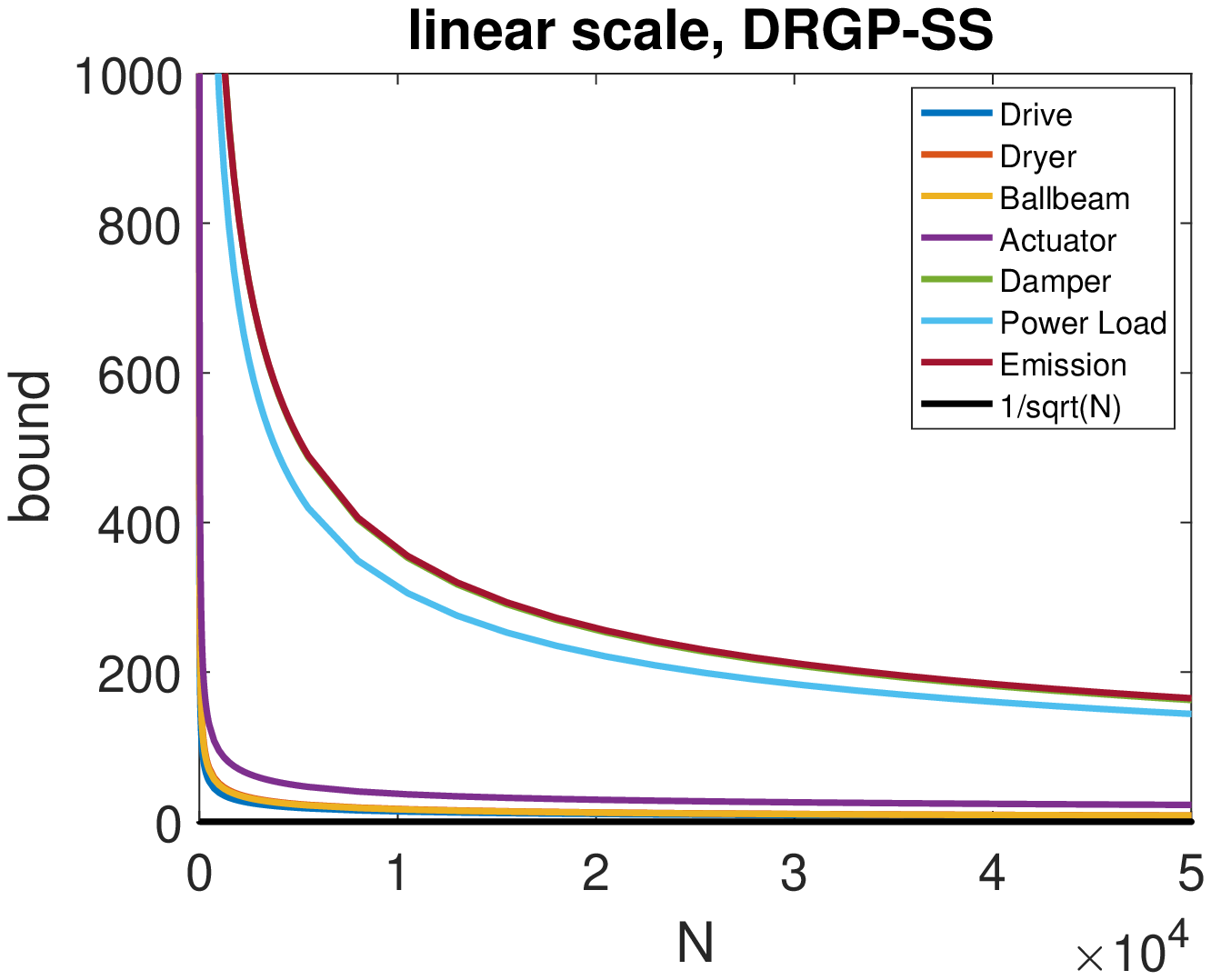}}\hspace{1cm}
  \includegraphics[width=0.4\textwidth]{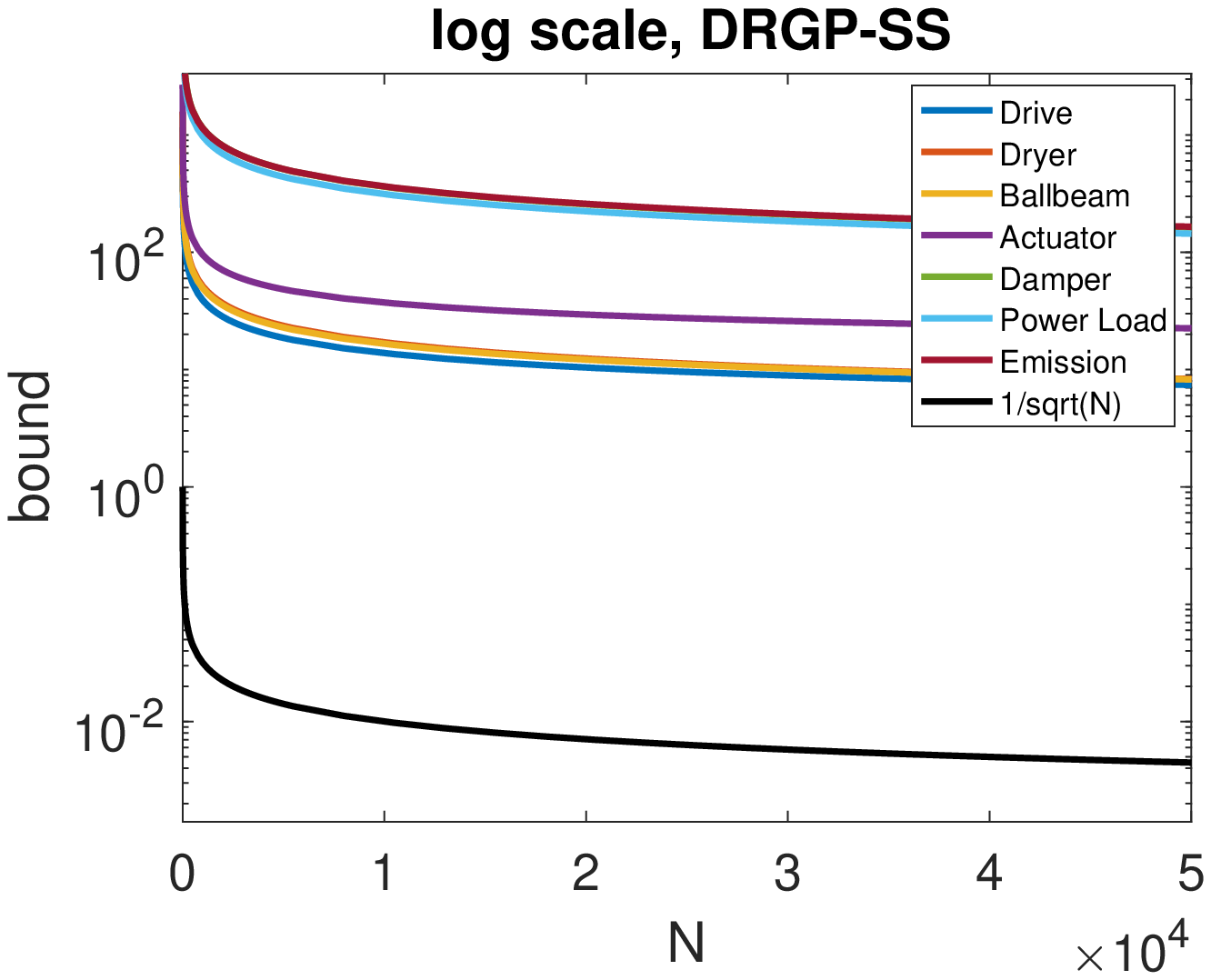}\\
	\includegraphics[width=0.4\textwidth]{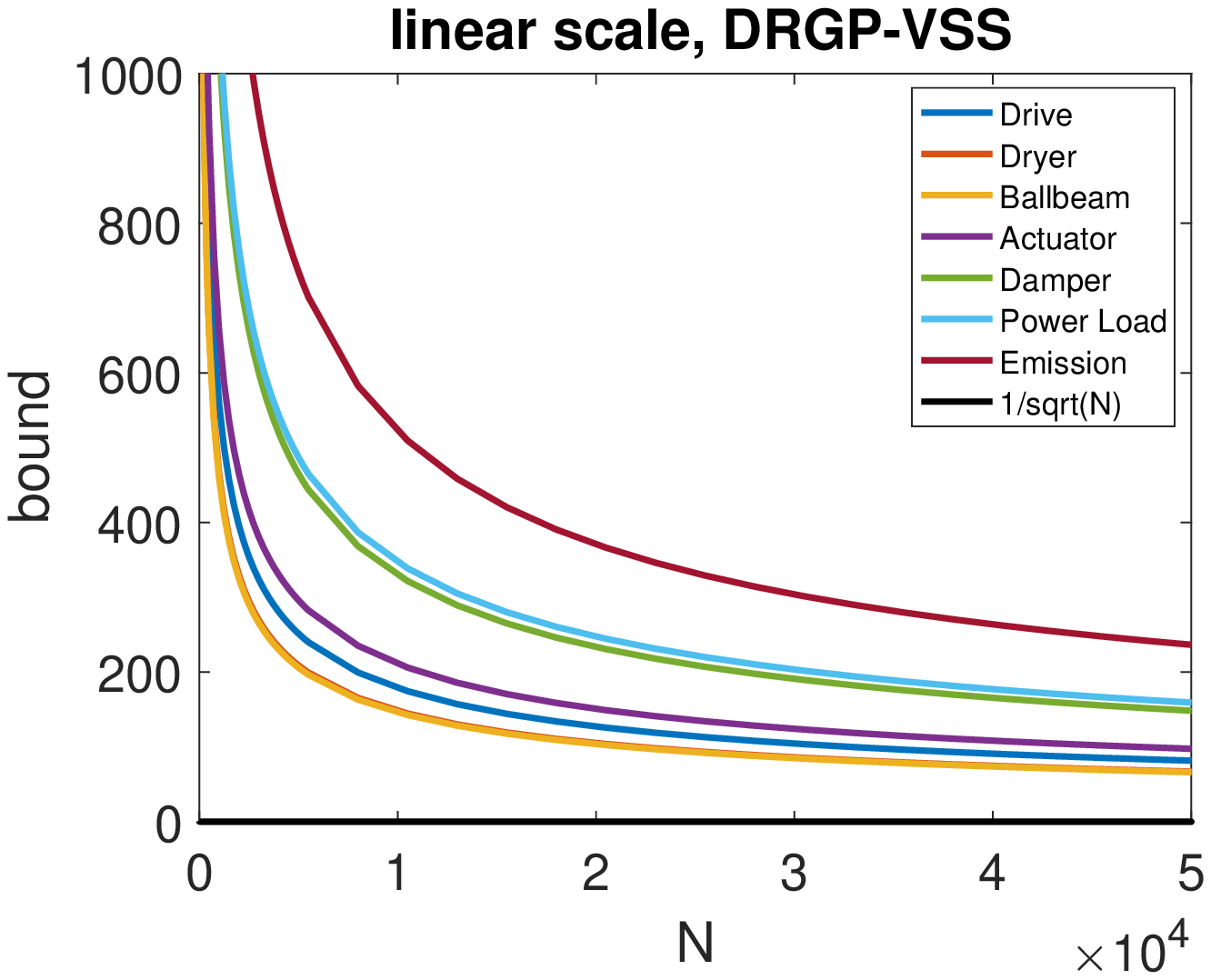}\hspace{1cm}
	\includegraphics[width=0.4\textwidth]{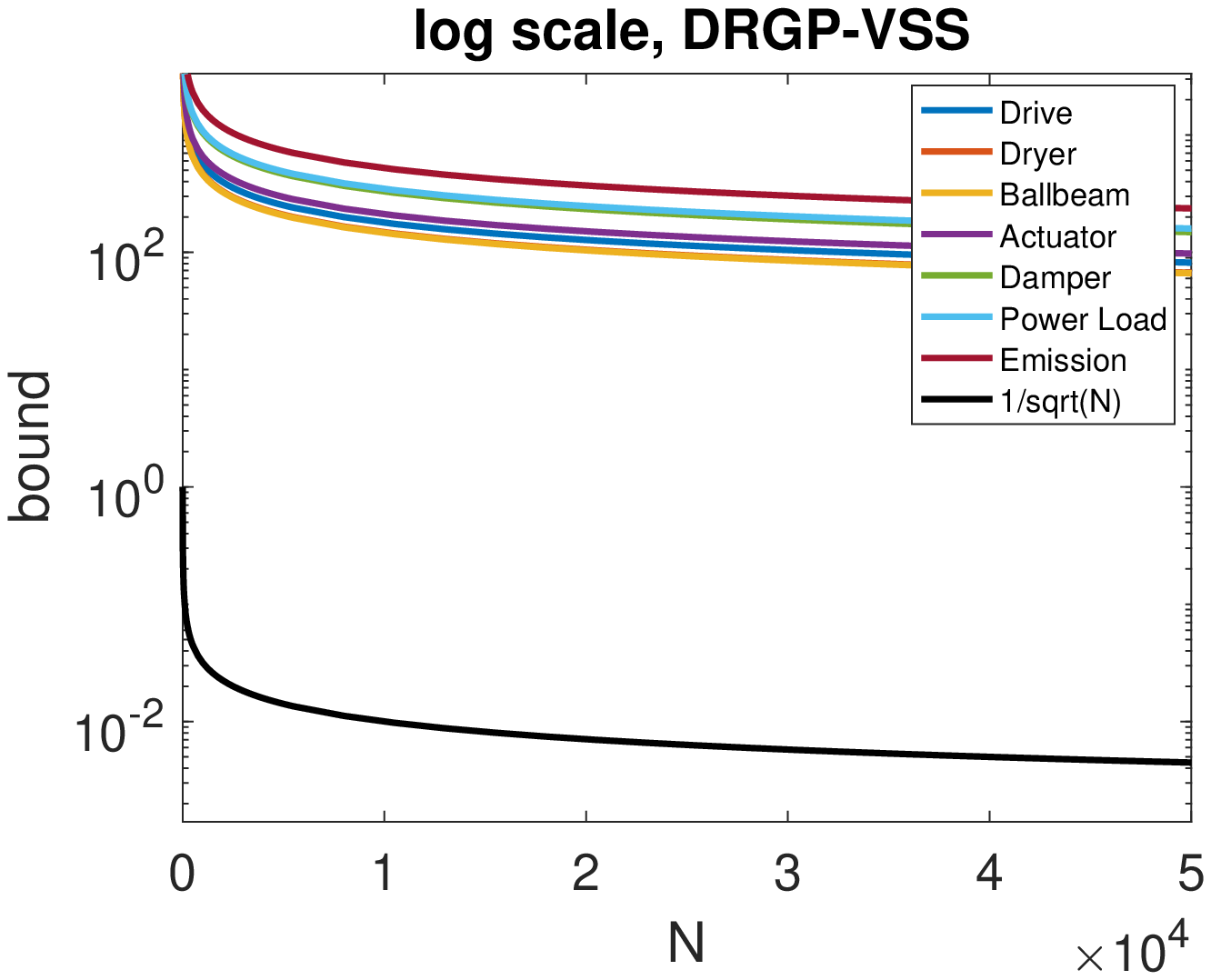}\\
	\includegraphics[width=0.4\textwidth]{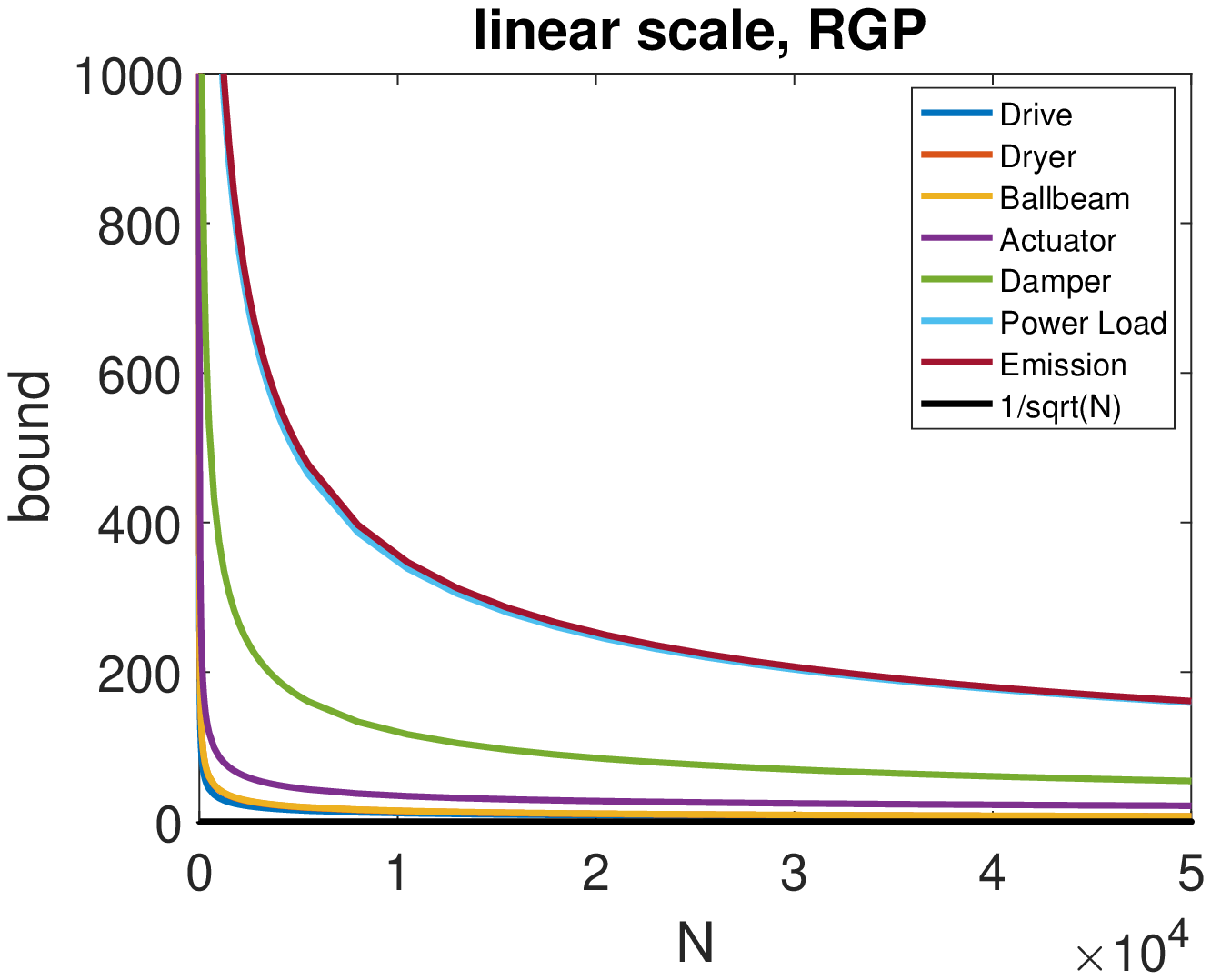}\hspace{1cm}
	\includegraphics[width=0.4\textwidth]{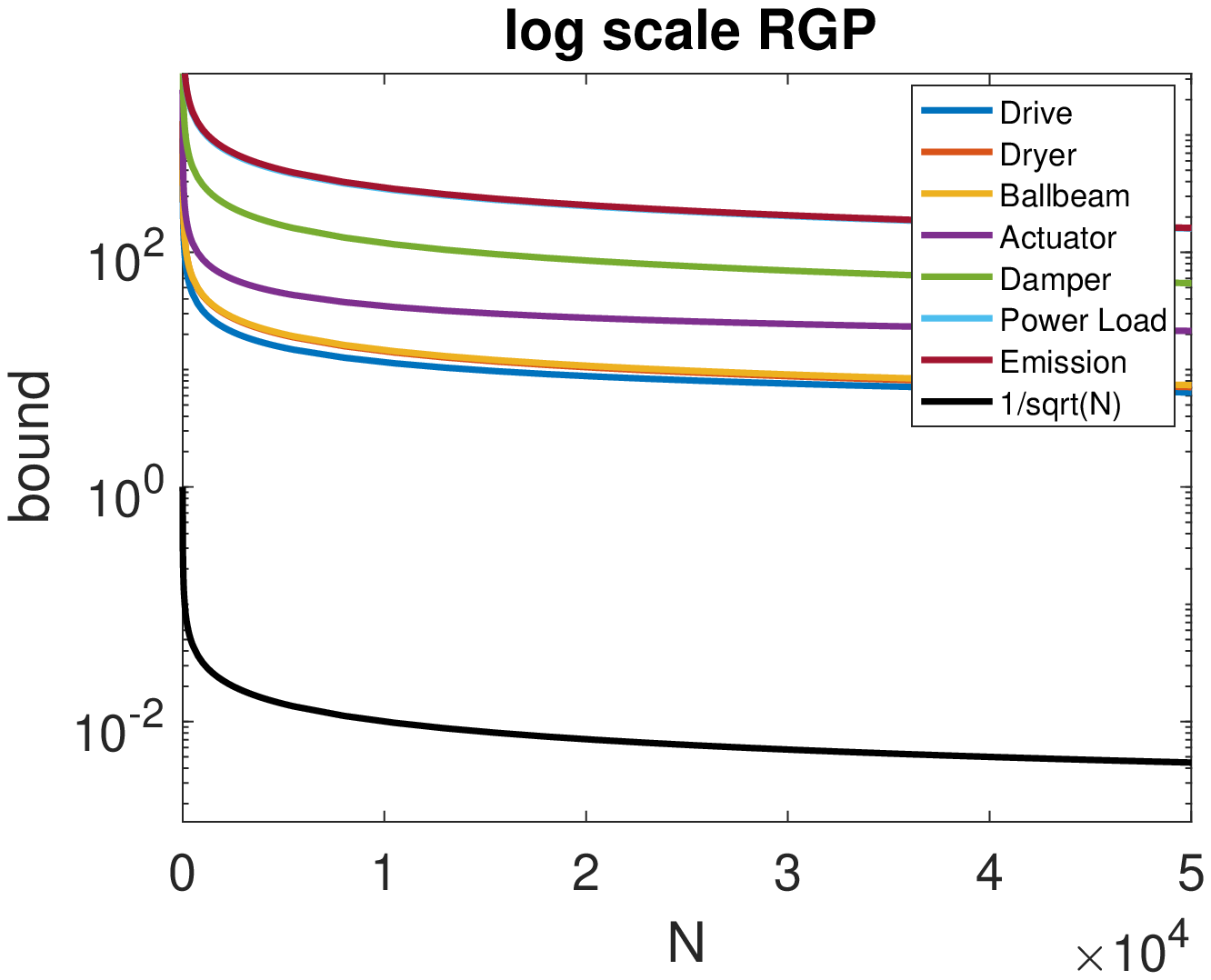}\\
	\caption{Bayesian-PAC bound evolution with $\lambda=\sqrt{N}$ for $N$ going from 1 to $50000$. Comparison for different data-sets and the DGPs called DRGP-(V)SS, RGP.}
\label{fig:VIS3}
\end{figure}
\noindent This experimental setup of creating from the existing training data-sets quasi-real data hints on generating ideal measured data in practice, because noise for data-based modeling can only be assumed to exist on a specific state, if one measures this state several times. For $Q_{\lambda}$ and the parameters without prior assumption we choose the already optimized parameters as in the experiments in \cite{foell2019deep}. In Figure \ref{fig:VIS2} in the Appendix~\ref{experimentsA} we see the test data-set Actuator and its prediction plus-minus two times standard deviation (SD), as well as its predicted hidden output-data plus-minus two times the SD. The covariance is assumed to be bounded with the model training result of $\textbf{Cov}[\boldsymbol{h}^{(l)}]$, the Lipschitz constant $\mathbf{S}^{(l)}=\max_{\mathsf{k}=1,...,\mathsf{K}-1}\bigl(|\mathrm{h}_{\mathsf{k}+1}^{\text{pred}}-\mathrm{h}_{\mathsf{k}}^{\text{pred}}|\bigr)$, and $\updelta^{(l)}=\dim(\mathbf{x}_\mathsf{k}^{(l)})$ for the layers $l=1,\dots,L+1$. In Figure \ref{fig:VIS3} we can now see the evolution of the bounds for the different models. We see a linear and a log scale with 7 different quasi-real data-sets created as described above. As expected the evolution is of order $\mathcal{O}\left(\frac{1}{\sqrt{N}}\right)$.
\begin{figure}[H]
\centering
	\includegraphics[width=0.4\textwidth]{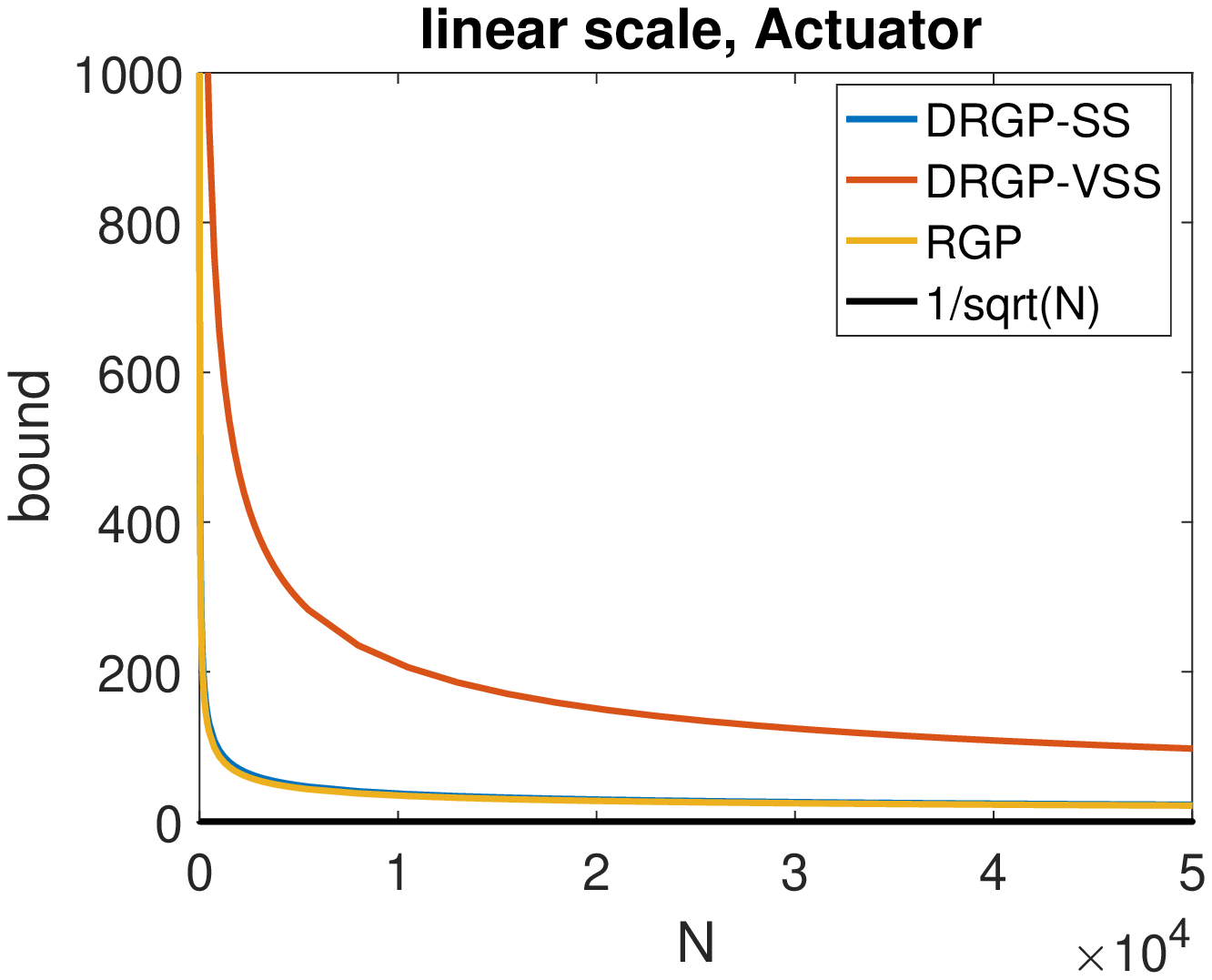}\hspace{1cm}
  \includegraphics[width=0.4\textwidth]{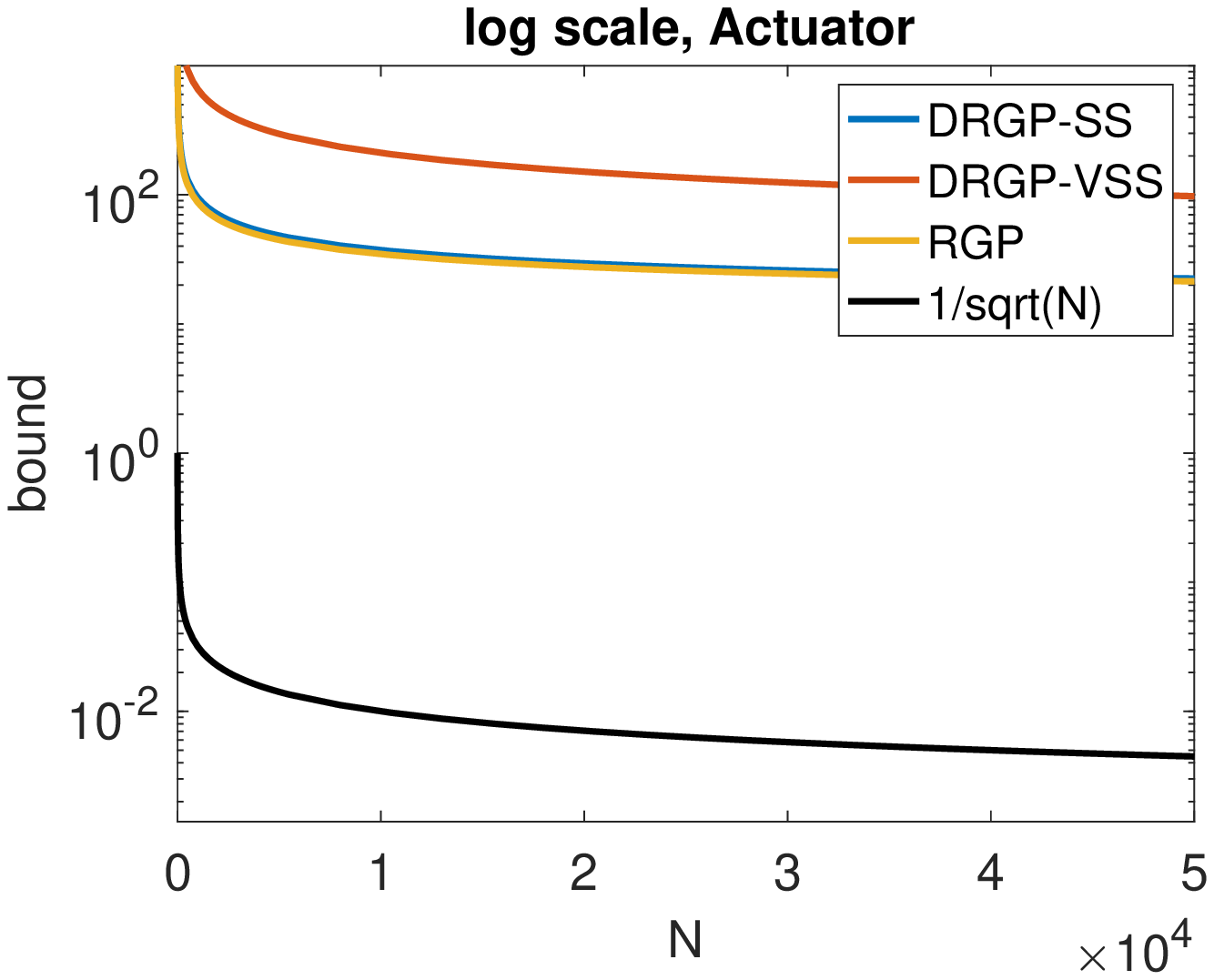}\\
	\includegraphics[width=0.4\textwidth]{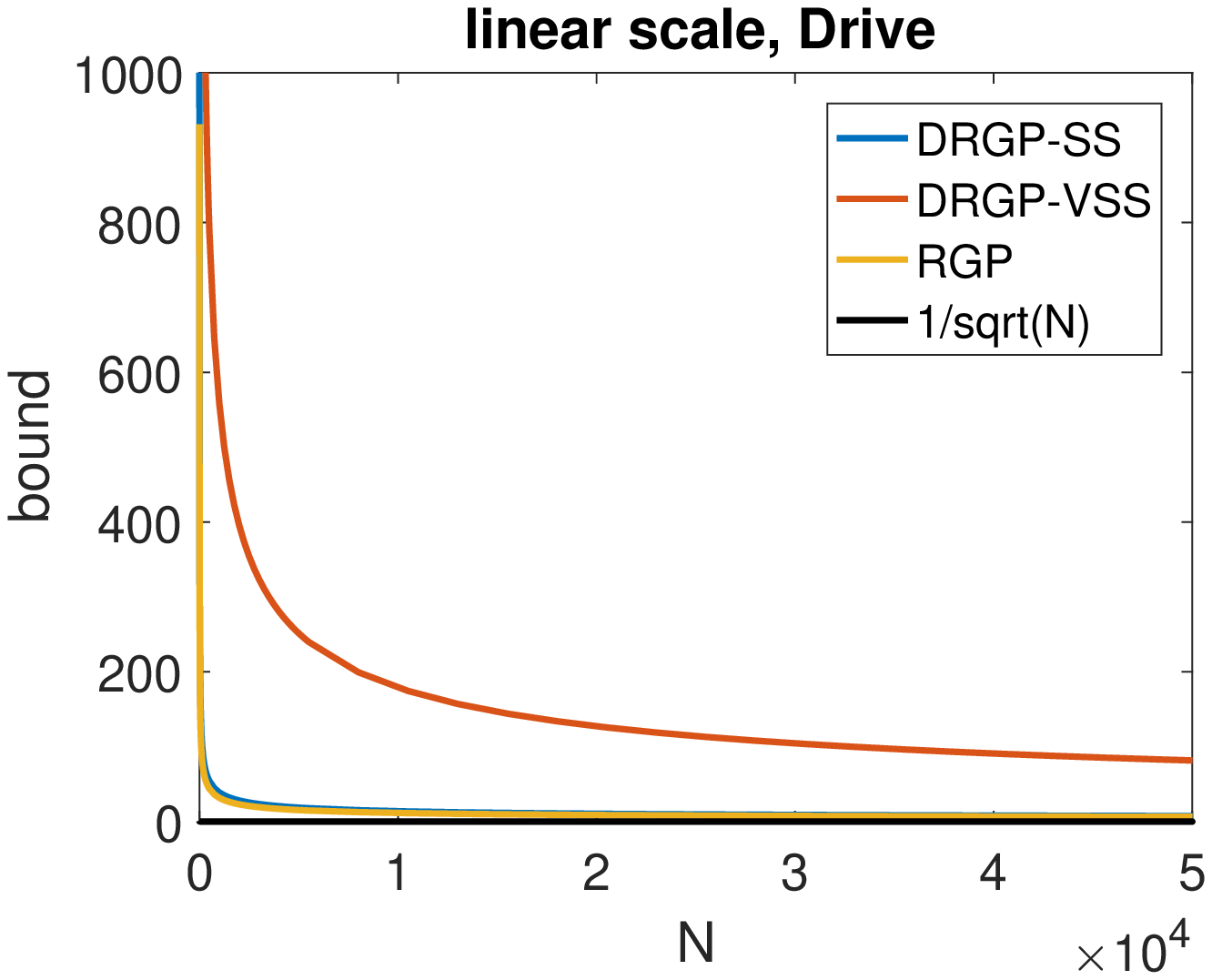}\hspace{1cm}
	\includegraphics[width=0.4\textwidth]{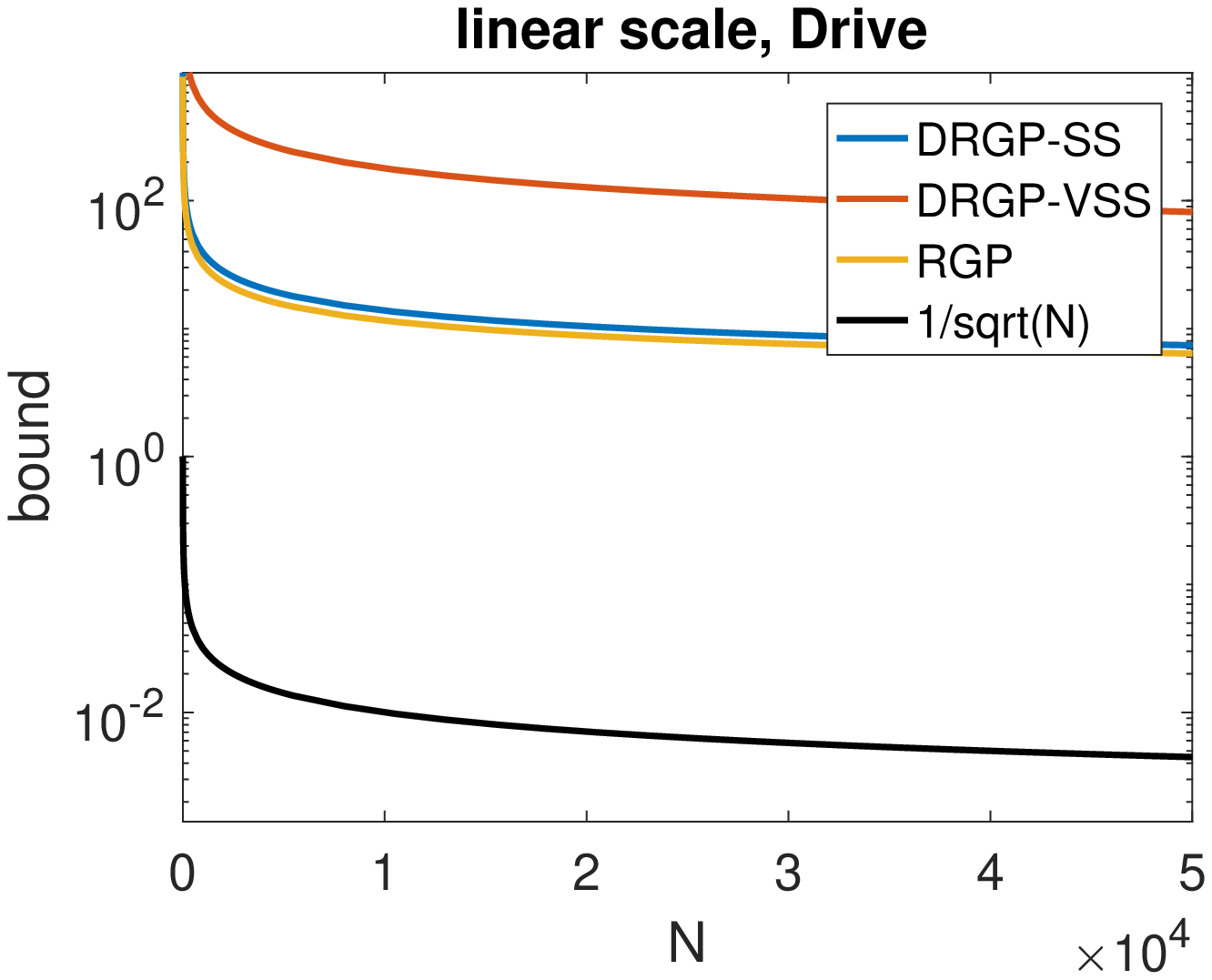}\\
\caption{Bayesian-PAC bound evolution with $\lambda=\sqrt{N}$ for $N$ going from 1 to $50000$. Comparison for different data-sets and the DGPs called DRGP-(V)SS, RGP.}
\label{fig:VIS4}	
\end{figure}
\noindent All three models reproduce similar curves, as seen Figure \ref{fig:VIS4}, which is not surprising, as we are dealing with trained parameters (the numerically minimized variational bound is present in the PAC-Bayesian bound, see Theorem 2 and Equation (\ref{PACVAR})). RGP, DRGP-SS compared to DRGP-VSS are more different resulting from fixing some parameters for DRGP-VSS during training.
\section{Conclusion}
\label{sec:conclusion}

The contribution of this paper is the derivation of PAC-Bayesian statements for several DGP models. Precise answers in different situations when modeling with DGPs e.g. designing controllers, are helpful in the overall modeling process. With the shown PAC-Bayesian statements we can give these answers. We see, that we can have the same intuition about the link of PAC-Bayesian bound and the marginal likelihood for our models, which use a variational bound instead of the the marginal likelihood. To show consistency, we extended the loss property to a quadratic-form-Gaussian loss function property, with convergence of order $\mathcal{O}\left(\frac{1}{\sqrt{N}}\right)$ for $\lambda=\sqrt{N}\to\infty$. We further showed several new forms of the bound in Theorem 5 and following, which give insights in the theory and practice. Our experiments show the evolution of the convergence for many data-sets used throughout the community for the DGPs of \cite{mattos2015recurrent} and \cite{foell2019deep} for the case of dynamic modeling.

\acks{We would like to acknowledge support for this project
from the ETAS GmbH.}


\vskip 0.2in
\bibliography{example_paper}
\newpage
\noindent{\large{\textbf{Appendix}}}\\

\noindent In this additional material we show some figures to the experimental results of Section~\eqref{experiments} and introduce the Definition of consistency and the union bound in Appendix \ref{definition} and the Proof of the Theorem 2, Theorem 3 in Appendix \ref{pacbayesian}, Theorem 5 in Appendix \ref{pacbayesian2} and the extension at the end of the paper in Appendix \ref{extensions}.
\appendix
\section{Experimental results}
\label{experimentsA}
\begin{figure}[H]
\centering
  \includegraphics[width=0.9\textwidth]{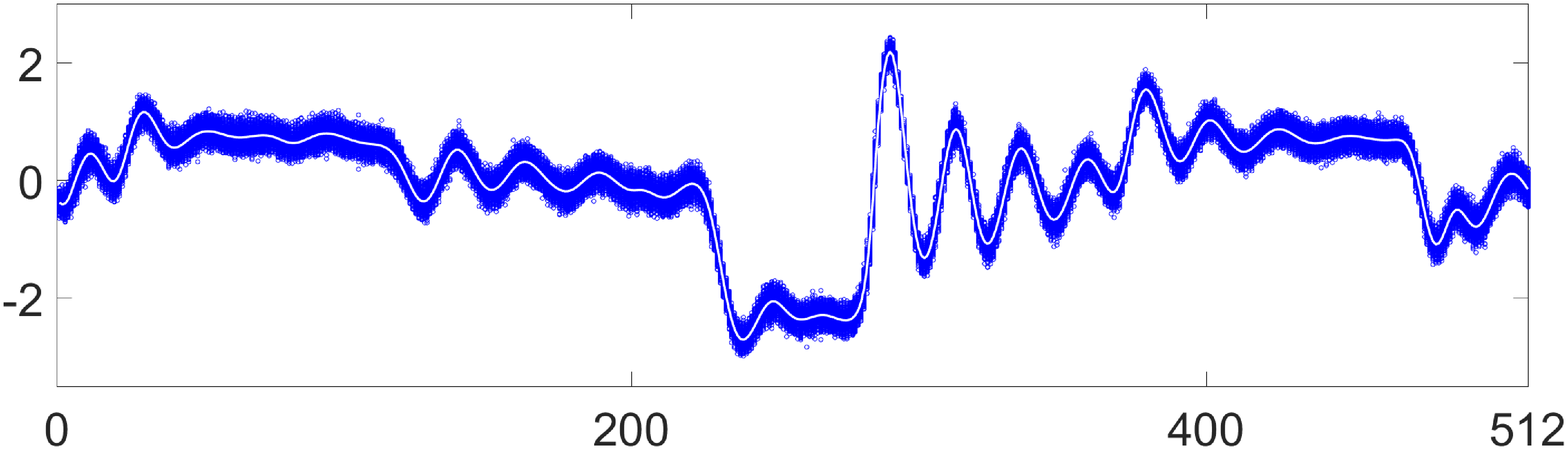}\\
	\caption{Quasi-real data for the training data-set Actuator, normalized data. Generation of the $N=50000$ samples with noise coming from the predictive posterior distribution after the training with $\mathsf{K}=512$ states. White: predicted training output data $\mathbf{y}^{\text{pred}}$, Blue: generated samples $\mathbf{y}^i$.}
\label{fig:VIS1}
	\includegraphics[width=0.4\textwidth]{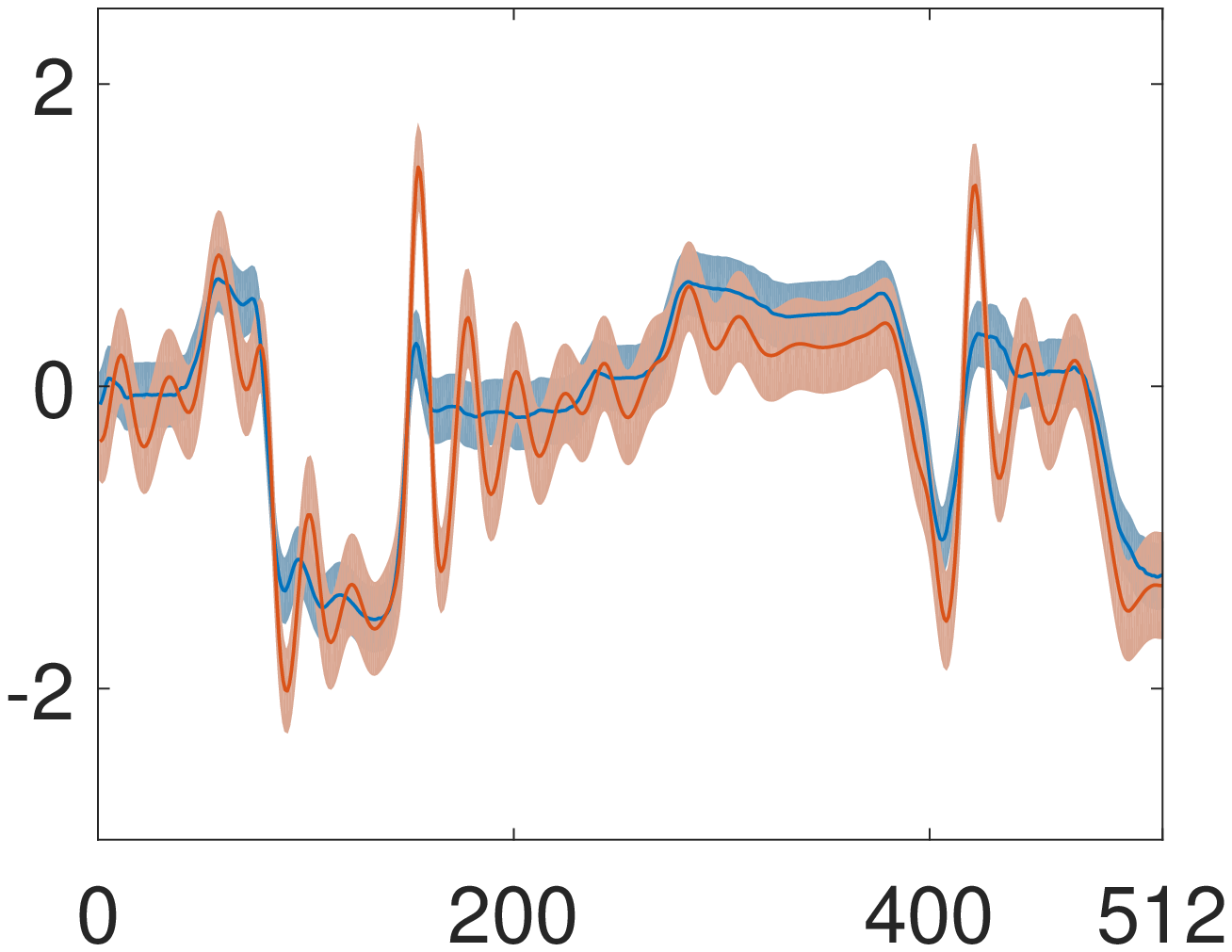}
	\includegraphics[width=0.4\textwidth]{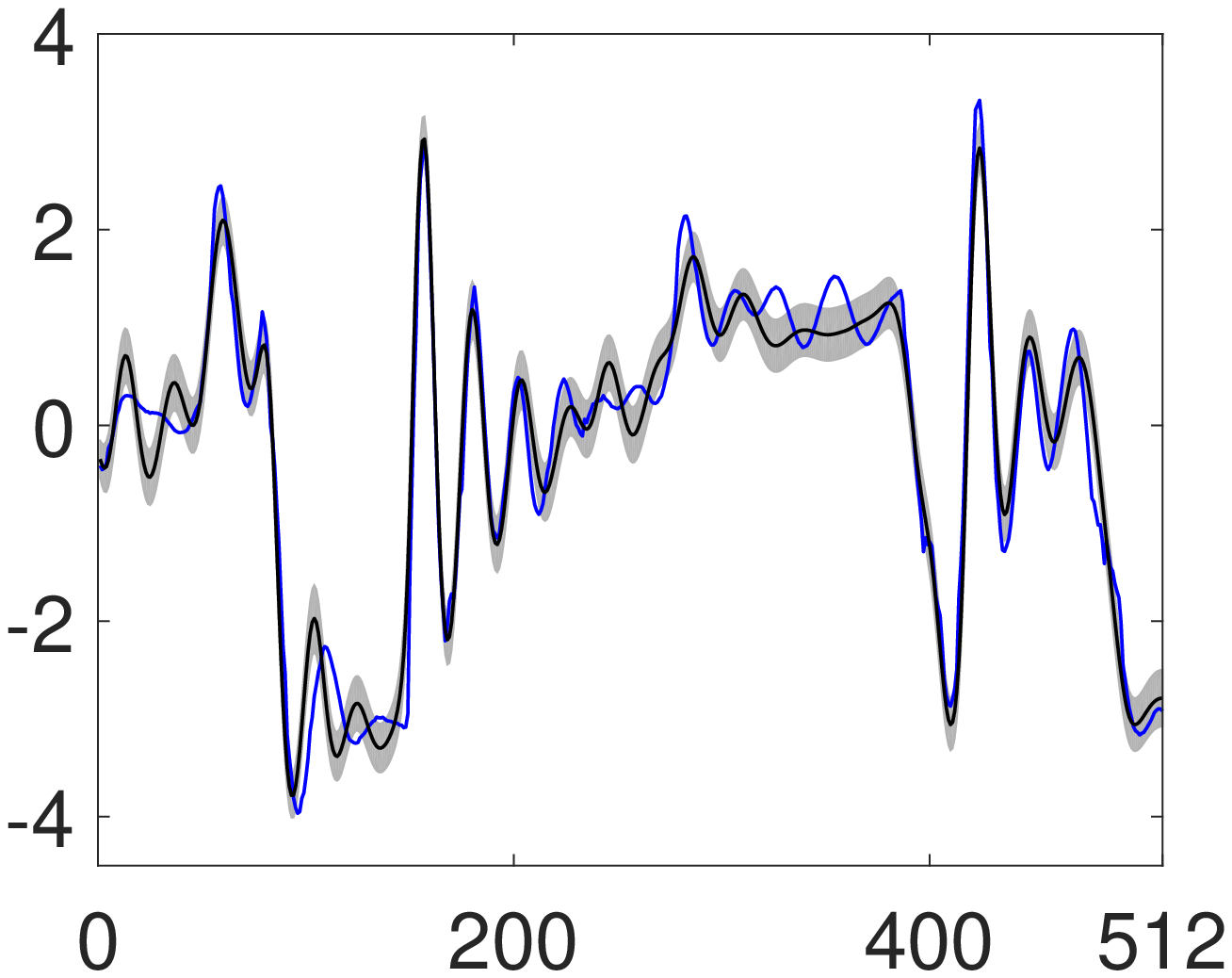}\\
	\caption{The Actuator data-set consists of $512$ training data and $512$ test data. Predicted test output-data for 2 hidden GP-layer on the left, Blue: 1. GP-layer $\pm 2$ times SD, Red: 2. GP-layer $\pm 2$ times SD. Predicted and real test output-data for output GP-layer on the right, Blue: real data, Black: prediction $\pm 2$  times SD. }
\label{fig:VIS2}	
\end{figure}
\newpage
\section{Definition Consisteny, Lemma 1. Union Bound and the two-side case}
\label{definition}
\begin{defi}[consistency]
We will say that a procedure that returns $\mathfrak{f}_{\uptheta}$ is consistent for a given measure $P_\mathsf{K}$, variational distribution $Q_{\text{\scalebox{.8}{\tiny{PAC}}}}$ and loss function $\ell$ if
\begin{align*} 
&\underset{{\hat\theta\sim Q_{\text{\scalebox{.8}{\tiny{PAC}}}}}}{\mathbf{E}}\left[\mathcal{L}_{ P_\mathsf{K}}^{\ell}(\mathfrak{f}_{\uptheta})\right]-\underset{{\hat\theta\sim Q_{\text{\scalebox{.8}{\tiny{PAC}}}}}}{\mathbf{E}}\left[\mathcal{L}_{\mathbb{D}}^{\ell}(\mathfrak{f}_{\uptheta})\right]\to 0,\; N\to\infty
\end{align*}
where convergence is assessed in a suitable manner, here in probability. If $\mathfrak{f}_{\uptheta}$ is consistent for all Borel probability measures then it is said to be universally consistent.
\end{defi}
\noindent{\bf Lemma 1. (Union Bound)} {\it Let $(\Sigma,\mathcal A, P)$ be a probability space and $\zeta_i$, $i=1,\dots,n$ be a set of countable events. We then have
\begin{align*} 
P\left(\bigcap_{i=1,\dots,n}\zeta_i\right)\leq \sum_{i=1}^n P\left(\zeta_i\right) - (n-1).
\end{align*}}
\\
\\
The statements of the \textit{two-side versions} of our statements follows by considering for the one-side statements the left hand-side of the inequality as two events $\zeta$, $-\zeta$ and by making the estimation for $\zeta\cap-\zeta$.\newpage
\section{Proof of Theorem 2 and Theorem 3}
\label{pacbayesian}
{\bf Proof}\\
We first introduce the model DRGP-(V)SS for the first page:\\
The priors $p_{\boldsymbol{a}^{(l)}}$, $p_{\boldsymbol{z}^{(l)}_m}$, $p_{\tilde{\boldsymbol{h}}^{(l)}}$ are 
\begin{align*}
&\boldsymbol{a}^{(l)}\sim\mathcal{N}(\mathbf{0},I_{M}),\quad\boldsymbol{z}_m^{(l)}\sim \mathcal{N}(\mathbf{0},I_Q),\quad\tilde{\boldsymbol{h}}^{(l)}\sim\mathcal{N}(\mathbf{0},I_{2H_{\mathrm{h}}-H_{\mathbf{x}}}),
\end{align*}
the product of them is defined as $P_{\text{\scalebox{.8}{\tiny{REV}}}}$, and the variational distributions $q_{\boldsymbol{a}^{(l)}}$, $q_{\boldsymbol{z}^{(l)}_m}$, $q_{h^{(l)}_\mathsf{k}}$ are
\begin{align*}
& \boldsymbol{a}^{(l)}\sim\mathcal{N}(\mathbf{m}^{(l)},\mathbf{s}^{(l)}),\quad\boldsymbol{z}_m^{(l)} \sim\mathcal{N}(\mathbf{\upalpha}_m^{(l)},\mathbf{\upbeta}_m^{(l)}),\quad h_\mathsf{k}^{(l)}\sim\mathcal{N}(\upmu^{(l)}_\mathsf{k},\uplambda^{(l)}_\mathsf{k}),
\end{align*}
the product of them is defined as $Q_{\text{\scalebox{.8}{\tiny{REV}}}}$, and where $\mathbf{\upbeta}_m^{(l)}\in \;\mathbb R^{Q\times Q}\;\text{is diagonal}$, for $\mathsf{k}=1,\dots,\mathsf{K}$, $m=1\dots,M$, $l=1,\dots,L+1$.\\ 
The upcoming statistics are $(\psi_1^{(l)})^T = \mathbf{E}_{q_{\boldsymbol{Z}^{(l)}}q_{h^{(l)}}}\left[\phi^{(l)}\right]\in\mathbb R^{M}$, $\Psi_1^{(l)} = \mathbf{E}_{q_{\boldsymbol{Z}^{(l)}}q_{\boldsymbol{h}^{(l)}}}\left[\Phi^{(l)}\right]\in\mathbb R^{\mathsf{K}\times M}$ and $\psi_2^{(l)} = \mathbf{E}_{q_{\boldsymbol{Z}^{(l)}}q_{h^{(l)}}}\left[\phi^{(l)}(\phi^{(l)})^T\right]\in\mathbb R^{M\times M}$, $\Psi_2^{(l)} = \mathbf{E}_{q_{\boldsymbol{Z}^{(l)}}q_{\boldsymbol{h}^{(l)}}}\left[(\Phi^{(l)})^T\Phi^{(l)}\right]\in\mathbb R^{M\times M}$, where with $\mathfrak{L} \stackrel{\mathrm{def}}= \text{diag}([2\pi \mathrm{l}_q]_{q=1}^Q)$, $\mathbf{p} \stackrel{\mathrm{def}}= [\mathrm{p}_1^{-1},\dots,\mathrm{p}_Q^{-1}]^T\in\mathbb R^{Q}$ and $\mathbf{Z}^{(l)} \stackrel{\mathrm{def}}=[\mathbf{z}_1^{(l)},\dots,\mathbf{z}_M^{(l)}]^T\in\mathbb R^{M\times Q}$ we have
\begin{align*}
\phi^{(l)}(\mathbf{x}^{(l)},\mathbf{Z}^{(l)}) \stackrel{\mathrm{def}} = \sqrt{2(\upsigma_{\text{power}}^{(l)})^2{M}^{-1}}\left[\cos(2\pi((\mathfrak{L}^{(l)})^{-1}\mathbf{z}_1^{(l)}+\mathbf{p}^{(l)})^T(\mathbf{x}^{(l)}-\mathbf{u}_1^{(l)}) + \mathrm{b}_1^{(l)}),\dots,\right.\\
\left.\cos(2\pi((\mathfrak{L}^{(l)})^{-1}\mathbf{z}_M^{(l)}+\mathbf{p}^{(l)})^T(\mathbf{x}^{(l)}-\mathbf{u}_M^{(l)}) + \mathrm{b}_M^{(l)})\right]^T\in\mathbb R^{M},
\end{align*}
sampling $b^{(l)}\sim\text{Unif}\left[0,2\pi\right]$ and $\boldsymbol{z}^{(l)}\sim\mathcal{N}(\mathbf{0},I_Q)$, where $K_{\mathsf{K}\mathsf{K}}^{(l),(\text{\tiny{SM}})} \stackrel{\mathrm{def}}= \Phi^{(l)}(\Phi^{(l)})^T$ and where $\Phi^{(l)} \stackrel{\mathrm{def}}= \left[\phi^{(l)}(\mathbf{x}_1^{(l)},\mathbf{Z}^{(l)}),\dots,\phi(\mathbf{x}_\mathsf{K}^{(l)},\mathbf{Z}^{(l)})\right]^T\in\mathbb R^{\mathsf{K}\times M}$.\\
\noindent We assume \textbf{Assumption 4} \textbf{(A4)}, that there exists parameters indexed with * (not necessarily unique) such that\\
\\
$\mathrm{y}^i = \psi^{(L+1)}_{1,*}\mathbf{m}^{(L+1)}_*+\epsilon^{y}_{i,*}$,\quad\text{and} $\mathrm{h}^{
i,{(l)}}=\mathrm{h}^{*,i,(l)}=\psi^{(l)}_{1,*}\mathbf{m}^{(l)}_{*}+\epsilon^{h^{(l)}}_{i,*}$,$\quad\text{for }\mathsf{K}=1$\\
$\mathbf{y}^i = \Psi^{(L+1)}_{1,*}\mathbf{m}^{(L+1)}_*+\epsilon^{\boldsymbol{y}}_{i,*}$,\quad\text{and} $\mathbf{h}^{
i,{(l)}}=\mathbf{h}^{*,i,(l)}=\Psi^{(l)}_{1,*}\mathbf{m}^{(l)}_{*}+\epsilon^{\mathbf{h}^{(l)}}_{i,*}$,$\quad\text{for }\mathsf{K}>1$\\
\\
for $i = 1,\dots,\bar{N}$, $\mathsf{k}=1,\dots,\mathsf{K},\quad l=1,\dots,L$, where the predictive posterior distributions of our DRGP-(V)SS models are\\
\\
$y_\mathsf{k}|\boldsymbol{m}_*^{(L+1)},{\boldsymbol\upalpha}^{(L+1)}_*,{\boldsymbol\upbeta}^{(L+1)}_*,{\boldsymbol\uplambda}^{(L)}_{\mathsf{k},*},\boldsymbol\upmu^{(L)}_{\mathsf{k},*}\sim\\
\mathcal{N}(\psi^{(L+1)}_{1,*,\mathsf{k}}\mathbf{m}^{(L+1)}_*,(\sigma_{\text{noise}*}^{(L+1)})^2+(\mathbf{m}^{(L+1)}_*)^T((\psi^{(L+1)}_{1,*,\mathsf{k}})^T\psi^{(L+1)}_{1,*,\mathsf{k}}-\psi^{(L+1)}_{2,*,\mathsf{k}})\mathbf{m}^{(L+1)}_*+\text{tr}(\psi^{(L+1)}_{2,*,\mathsf{k}}\mathbf{s}^{(L+1)}_*))$\\
\\
$h^{(l)}_\mathsf{k}|\boldsymbol{m}_*^{(l)},{\boldsymbol\upalpha}^{(l)}_*,{\boldsymbol\upbeta}^{(l)}_*,{\boldsymbol\uplambda}^{(l)}_{\mathsf{k},*},\boldsymbol\upmu^{(l)}_{\mathsf{k},*},{\boldsymbol\uplambda}^{(l-1)}_{\mathsf{k},*},\boldsymbol\upmu^{(l-1)}_{\mathsf{k},*}\sim\\
\mathcal{N}(\psi^{(l)}_{1,*,\mathsf{k}}\mathbf{m}^{(l)}_*,(\sigma_{\text{noise}*}^{(l)})^2+(\mathbf{m}^{(l)}_*)^T((\psi^{(l)}_{1,*,\mathsf{k}})^T\psi^{(l)}_{1,*,\mathsf{k}}-\psi^{(l)}_{2,*,\mathsf{k}})\mathbf{m}^{(l)}_*+\text{tr}(\psi^{(l)}_{2,*,\mathsf{k}}\mathbf{s}^{(l)}_*),\quad l=2,\dots,L$,\\
\\
$h^{(1)}_\mathsf{k}|\boldsymbol{m}_*^{(1)},{\boldsymbol\upalpha}^{(1)}_*,{\boldsymbol\upbeta}^{(l)}_*,{\boldsymbol\uplambda}^{(1)}_{\mathsf{k},*},\boldsymbol\upmu^{(1)}_{\mathsf{k},*},\mathbf{\bar{x}}_{\mathsf{k}}\sim\\
\mathcal{N}(\psi^{(1)}_{1,*,\mathsf{k}}\mathbf{m}^{(1)}_*,(\sigma_{\text{noise}*}^{(1)})^2+(\mathbf{m}^{(1)}_*)^T((\psi^{(1)}_{1,*,\mathsf{k}})^T\psi^{(1)}_{1,*,\mathsf{k}}-\psi^{(l)}_{2,*,\mathsf{k}})\mathbf{m}^{(1)}_*+\text{tr}(\psi^{(1)}_{2,*,\mathsf{k}}\mathbf{s}^{(1)}_*)$.\\
\\
Further we identify for the distribution $P_\mathsf{K}$ the following density\\
{\centering
$p_{\boldsymbol{y}|\boldsymbol{m}_*^{(L+1)},{\boldsymbol\upalpha}^{(L+1)}_*,{\boldsymbol\upbeta}^{(L+1)}_*,{\boldsymbol\uplambda}^{(L)}_*,\boldsymbol\upmu^{(L)}_*}p_{\boldsymbol{h}^{(1)}|\boldsymbol{m}_*^{(1)},{\boldsymbol\upalpha}^{(1)}_*,{\boldsymbol\upbeta}^{(l)}_*,{\boldsymbol\uplambda}^{(1)}_{\mathsf{k},*},\boldsymbol\upmu^{(1)}_{\mathsf{k},*},\mathbf{\bar{x}}_{\mathsf{k}}}\prod\nolimits_{l=2}^{L}p_{\boldsymbol{h}^{(l)}|\boldsymbol{m}_*^{(l)},{\boldsymbol\upalpha}^{(l)}_*,{\boldsymbol\upbeta}^{(l)}_*,{\boldsymbol\uplambda}^{(l)}_*,\boldsymbol\upmu^{(l)}_*,{\boldsymbol\uplambda}^{(l-1)}_*,\boldsymbol\upmu^{(l-1)}_*}$.}\\
\\
This assumption is a very reasonable assumption, because otherwise we would use another model for the modeling problem. One has to take care of the fact, that these parameters do not necessarily need to be equal to the ones produced by the Bayes risk minimization $\mathcal{L}_{ P_\mathsf{K}}^{\ell_{\text{nll}}}(\mathfrak{f}_{\scalebox{.8}{$\bar{\uptheta}$}})=\text{arg}\min\limits_{\hat\uptheta\in\hat{\Theta}}\left(\mathcal{L}_{P_\mathsf{K}}^{\ell_{\text{nll}}}(\mathfrak{f}_{\scalebox{.8}{$\uptheta$}})\right)$.\\
For simplicity we did not show here the optimal parameters $\sigma_{*\text{power}}^{(l)},\frac{diag(\mathfrak{L}^{(l)}_{*})}{2\pi}$, $\mathbf{U}^{(l)}_{*}=[\mathbf{u}_{1,*}^{(l)},\dots,\mathbf{u}_{M,*}^{(l)}]^T$, $\mathbf{p}^{(l)}_{*}$, $\mathbf{b}^{(l)}_{*}=[\mathrm{b}_{1,*}^{(l)},\dots,\mathrm{b}_{M,*}^{(l)}]^T$ without prior assumption. We further assume that we use for all models the same sparsity parameter $M$ and for the parameters star the optimal distribution $\mathbf{m}^{(L+1)}_{*}=(A^{(L+1)})^{-1}\left(\Psi_{1,*}^{(L+1)}\right)^T\mathbf{Y}$ and $\mathbf{m}^{(l)}_{*}=(A^{(l)})^{-1}\left(\Psi_{1,*}^{(l)}\right)^T\boldsymbol\upmu^{(l)}_{*}$ for $l=1,\dots,L$, as well as $\mathbf{s}^{(l)}_{*}=(\sigma_{\text{noise}*}^{(l)})^2(A^{(l)})^{-1}$ for $A^{(l)} = (\Psi_{2,*}^{(l)} + (\sigma_{\text{noise}*}^{(l)})^2)$ for $l=1,\dots,L+1$.\\
The theorems are adaptable for~\cite{mattos2015recurrent} by using the definition of the statistics $\Psi$ from their paper, as well as the definitions of their priors, variational distributions and $\mathbf{m}^{(L+1)}_{*}=(B^{(L+1)})^{-1}\left(\Psi_{1,*}^{(L+1)}\right)^T\mathbf{Y}$, $\mathbf{m}^{(l)}_{*}=(B^{(l)})^{-1}\left(\Psi_{1,*}^{(l)}\right)^T\boldsymbol\upmu^{(l)}_{*}$ for $l=1,\dots,L$, and $\mathbf{s}^{(l)}_{*}=-\Psi_{0,*}^{(l)} + (K_{MM*}^{(l)})^{-1}-(\sigma_{\text{noise}*}^{(l)})^2(B^{(l)})^{-1}$ for $B^{(l)} = (\Psi_{2,*}^{(l)} + K_{MM*}^{(l)}(\sigma_{\text{noise}*}^{(l)})^2)$ for $l=1,\dots,L+1$.\\
\\
\textit{\textbf{Discussion}} \textbf{(D1)}: \textit{Filling in here the optimal distribution for~\cite{foell2019deep} and~\cite{mattos2015recurrent} after taking the expectation w.r.t $P_K$ needs a explanation, because the optimal distribution depends on $\mathbf{Y}$. As our models \cite{foell2019deep} and \cite{mattos2015recurrent} inherently make use of predicting the mean values and variance depending on the output data $\mathbf{Y}$, these values are assumed to be constant and not random variables. Theoretically this should be noticed as a special characteristic of such a model, which takes measurements to derive a dependent probability distribution as model.}\\
\\\newpage
\noindent We proof the bound and loss property for $\Psi^{\ell}(\lambda,N)$ with prior $P_{\text{\scalebox{.8}{\tiny{PAC}}}}$ equal $P_{\text{\scalebox{.8}{\tiny{REV}}}}$ and $Q_{\text{\scalebox{.8}{\tiny{PAC}}}}$ equal variational distribution $Q_{\text{\scalebox{.8}{\tiny{REV}}}}$ and $\mathbf{X}^{(1)}=[\mathbf{x}_1^{(1)},\dots,\mathbf{x}_\mathsf{K}^{(1)}]$ is set by an experimenter. First we show an abbreviation for $\mathcal{L}_{P_\mathsf{K}}^{\ell_{\text{nll}}}(\mathfrak{f}_{\uptheta})$ in $\Psi^{\ell}(\lambda,1)$:
\begin{flalign*} 
\mathcal{L}_{P_\mathsf{K}}^{\ell_{\text{nll}}}(\mathfrak{f}_{\uptheta})&=\underset{\boldsymbol{y}\sim{P_\mathsf{K}}}{\mathbf{E}}\left[\ell_{\text{nll}}(\mathfrak{f}_{\uptheta},\mathbf{y})\right]\\
& = \underset{\boldsymbol{y}\sim{P_\mathsf{K}}}{\mathbf{E}}\left[-\left(\log(p(\mathbf{y}|\uptheta^{(L+1)},\mathbf{X}^{(L+1)}))+\sum\nolimits_{l=1}^{L}\log(p(\mathbf{h}^{(l)}|\uptheta^{(l)},,\mathbf{X}^{(l)}))\right)\right]\\
& = \underset{\boldsymbol{y}\sim{P_\mathsf{K}}}{\mathbf{E}}\left[\frac{(\mathbf{y}- \Phi^{(L+1)}\mathbf{a}^{(L+1)})^T(\mathbf{y}- \Phi^{(L+1)}\mathbf{a}^{(L+1)})}{2(\sigma_{\text{noise}}^{(L+1)})^2}+\mathsf{K}\frac{\log\left(2\pi(\sigma_{\text{noise}}^{(l)})^2\right)}{2}\right.\numberthis\label{empirical}\\
&\left.+\sum\nolimits_{l=1}^L \frac{(\mathbf{h}^{(l)}- \Phi^{(l)}\mathbf{a}^{(l)})^T(\mathbf{h}^{(l)}- \Phi^{(l)}\mathbf{a}^{(l)})}{2(\sigma_{\text{noise}}^{(l)})^2}+\mathsf{K}\frac{\log\left(2\pi(\sigma_{\text{noise}}^{(l)})^2\right)}{2}\right].
\end{flalign*}
Then it could be replaced by the following equations
\begin{align*}
&=\frac{\underset{\boldsymbol{y}\sim{P_\mathsf{K}}}{\mathbf{E}}\left[\mathbf{y}^T\mathbf{y}-2\mathbf{y}^T\Psi^{(L+1)}_{1,*}\mathbf{m}^{(L+1)}_*+(\Psi^{(L+1)}_{1,*}\mathbf{m}^{(L+1)}_*)^T(\Psi^{(L+1)}_{1,*}\mathbf{m}^{(L+1)}_*)\right]}{2(\sigma_{\text{noise}}^{(L+1)})^2}\\
&+\frac{(\Psi^{(L+1)}_{1,*}\mathbf{m}^{(L+1)}_*)^T(\Psi^{(L+1)}_{1,*}\mathbf{m}^{(L+1)}_*)}{2(\sigma_{\text{noise}}^{(L+1)})^2}-\frac{\underset{\boldsymbol{y}\sim{P_\mathsf{K}}}{\mathbf{E}}\left[2\mathbf{y}^T\Phi^{(L+1)}\mathbf{a}^{(L+1)}\right]}{2(\sigma_{\text{noise}}^{(L+1)})^2}\\
&+\frac{(\Phi^{(L+1)}\mathbf{a}^{(L+1)})^T(\Phi^{(L+1)}\mathbf{a}^{(L+1)})}{2(\sigma_{\text{noise}}^{(L+1)})^2}+\mathsf{K}\frac{\log\left(2\pi(\sigma_{\text{noise}}^{(L+1)})^2\right)}{2}\\
&+\sum\nolimits_{l=1}^L \frac{\underset{\boldsymbol{y}\sim{P_\mathsf{K}}}{\mathbf{E}}\left[(\mathbf{h}^{(l)})^T\mathbf{h}^{(l)}-2(\mathbf{h}^{(l)})^T\Psi^{(l)}_{1,*}\mathbf{m}^{(l)}_*+(\Psi^{(l)}_{1,*}\mathbf{m}^{(l)}_*)^T\Psi^{(l)}_{1,*}\mathbf{m}^{(l)}_*\right]}{2(\sigma_{\text{noise}}^{(l)})^2}\\
&+\frac{(\Psi^{(l)}_{1,*}\mathbf{m}^{(l)}_*)^T\Psi^{(l)}_{1,*}\mathbf{m}^{(l)}_*}{2(\sigma_{\text{noise}}^{(l)})^2} -\frac{\underset{\boldsymbol{y}\sim{P_\mathsf{K}}}{\mathbf{E}}\left[2(\mathbf{h}^{(l)})^T\Phi^{(l)}\mathbf{a}^{(l)}\right]}{2(\sigma_{\text{noise}}^{(l)})^2}+\frac{\underset{\boldsymbol{y}\sim{P_\mathsf{K}}}{\mathbf{E}}\left[(\Phi^{(l)}\mathbf{a}^{(l)})^T\Phi^{(l)}\mathbf{a}^{(l)}\right]}{2(\sigma_{\text{noise}}^{(l)})^2}\\
&+\mathsf{K}\frac{\log\left(2\pi(\sigma_{\text{noise}}^{(l)})^2\right)}{2}\\
&=\sum\nolimits_{l=1}^{L+1}\frac{(\mathbf{m}^{(l)}_*)^T(\Psi^{(l)}_{2,*}-(\Psi^{(l)}_{1,*})^T\Psi^{(l)}_{1,*})\mathbf{m}^{(l)}_*}{2(\sigma_{\text{noise}}^{(l)})^2}+\text{tr}\left(\frac{\Psi^{(l)}_{2,*}\mathbf{s}^{(l)}_*}{2(\sigma_{\text{noise}}^{(l)})^2}\right)+\mathsf{K}\frac{(\sigma_{\text{noise}*}^{(l)})^2}{2(\sigma_{\text{noise}}^{(l)})^2}\\
&+\frac{(\Psi^{(l)}_{1,*}\mathbf{m}^{(l)}_*-\Phi^{(l)}\mathbf{a}^{(l)})^T(\Psi^{(l)}_{1,*}\mathbf{m}^{(l)}_*-\Phi^{(l)}\mathbf{a}^{(l)})}{2(\sigma_{\text{noise}}^{(l)})^2}+\mathsf{K}\frac{\log\left(2\pi(\sigma_{\text{noise}}^{(l)})^2\right)}{2}\\
\end{align*}
or just plugging in the \textbf{(A4)} for the samples $\mathrm{y}$, $\mathrm{h}^{(l)}$ and calculating the expectation.\\
\\
If we assume the functions\\
\begin{flalign*}
&\mathfrak{f}^{(1)}:\mathbb R^{Q_r^{(1)}}\to\mathbb R,\numberthis\label{predictor1}\\
&(\mathbf{a}^{(1)},\mathbf{z}^{(1)},\mathbf{x}^{(1)},M,\sigma_{\text{power},{\mathfrak{f}}}^{(1)},\frac{\text{diag}(\mathfrak{L}^{(1)}_\mathfrak{f})}{2\pi},\mathbf{p}^{(1)}_\mathfrak{f},\mathbf{u}^{(1)}_\mathfrak{f},\mathbf{b}^{(1)}_\mathfrak{f})^T\to\phi^{(1)}\mathbf{a}^{(1)}\\
&\mathfrak{f}^{(l)}:\mathbb R^{Q_r^{(l)}}\to\mathbb R,\numberthis\label{predictor2}\\
&(\mathbf{a}^{(l)},\mathbf{z}^{(l)},\mathbf{x}^{(l)},M,\sigma_{\text{power},{\mathfrak{f}}}^{(l)},\frac{\text{diag}(\mathfrak{L}^{(l)}_\mathfrak{f})}{2\pi},\mathbf{p}^{(l)}_\mathfrak{f},\mathbf{u}^{(l)}_\mathfrak{f},\mathbf{b}^{(l)}_\mathfrak{f})^T\to\phi^{(l)}\mathbf{a}^{(l)},\quad l=2,\dots,L\\
&\mathfrak{f}^{(L+1)}:\mathbb R^{Q_r^{(L+1)}}\to\mathbb R,\numberthis\label{predictor3}\\
&(\mathbf{a}^{(L+1)},\mathbf{z}^{(L+1)},\mathbf{x}^{(L+1)},M,\sigma_{\text{power},{\mathfrak{f}}}^{(L+1)},\frac{\text{diag}(\mathfrak{L}^{(L+1)}_\mathfrak{f})}{2\pi},\mathbf{p}^{(L+1)}_\mathfrak{f},\mathbf{u}^{(L+1)}_\mathfrak{f},\mathbf{b}^{(L+1)}_\mathfrak{f})^T\\
&\to\phi^{(L+1)}\mathbf{a}^{(L+1)}\\
&\mathfrak{v}^{(1)}_{\mathbf{\bar{x}}_{\mathsf{k}}}:\mathbb R^{Q_v^{(1)}}\to\mathbb R,\numberthis\label{predictor4}\\
&(\mathbf{m}^{(1)},{\boldsymbol\upalpha}^{(1)},{\boldsymbol\upbeta}^{(1)},\boldsymbol\upmu^{(1)}_{\mathsf{k}},\boldsymbol\uplambda^{(1)}_{\mathsf{k}},M,\sigma_{\text{power},{\mathfrak{v}}}^{(1)},\frac{\text{diag}(\mathfrak{L}^{(1)}_\mathfrak{v})}{2\pi},\mathbf{p}^{(1)}_\mathfrak{v},\mathbf{u}^{(1)}_\mathfrak{v},\mathbf{b}^{(1)}_\mathfrak{v})^T\to\psi^{(1)}_{1,{\mathsf{k}}}\mathbf{m}^{(1)}\\
&\mathfrak{v}^{(l)}:\mathbb R^{Q_v^{(l)}}\to\mathbb R,\numberthis\label{predictor5}\\
&(\mathbf{m}^{(l)},{\boldsymbol\upalpha}^{(l)},{\boldsymbol\upbeta}^{(l)},\boldsymbol\uplambda^{(l)}_{\mathsf{k}},\boldsymbol\upmu^{(l)}_{\mathsf{k}},\boldsymbol\uplambda^{(l-1)}_{\mathsf{k}},\boldsymbol\upmu^{(l-1)}_{\mathsf{k}},M,\sigma_{\text{power},{\mathfrak{v}}}^{(l)},\frac{\text{diag}(\mathfrak{L}^{(l)}_\mathfrak{v})}{2\pi},\mathbf{p}^{(l)}_\mathfrak{v},\mathbf{u}^{(l)}_\mathfrak{v},\mathbf{b}^{(l)}_\mathfrak{v})^T\\
&\to\psi^{(l)}_{1,{\mathsf{k}}}\mathbf{m}^{(l)},\quad l=2,\dots,L\\
&\mathfrak{v}^{(L+1)}:\mathbb R^{Q_v^{(L+1)}}\to\mathbb R,\numberthis\label{predictor5}\\
&(\mathbf{m}^{(L+1)},{\boldsymbol\upalpha}^{(L+1)},{\boldsymbol\upbeta}^{(L+1)},\boldsymbol\upmu^{(L)}_{\mathsf{k}},\boldsymbol\uplambda^{(L)}_{\mathsf{k}},M,\sigma_{\text{power},{\mathfrak{v}}}^{(L+1)},\frac{\text{diag}(\mathfrak{L}^{(L+1)}_\mathfrak{v})}{2\pi},\mathbf{p}^{(L+1)},\mathbf{u}^{(L+1)}_\mathfrak{v},\mathbf{b}^{(L+1)}_\mathfrak{v})^T\\
&\to\psi^{(L+1)}_{1,{\mathsf{k}}}\mathbf{m}^{(L+1)}
\end{flalign*}
to be stochastic Lipschitz \textbf{(A1)} with constants $\mathbf{S}^{(l)}$.\\
This means for two functions $\mathfrak{v}^{(l)}:\mathbb R^{p+v}\to\mathbb R$, $\mathfrak{f}^{(l)}:\mathbb R^p\to\mathbb R$, where one is a stochastic version of the other (some parameters in the function $\mathfrak{f}^{(l)}$ are assumed to have a mean $\uptheta_\mathrm{m}^{(l)}$ and variance $\uptheta_\mathrm{v}^{(l)}$, which results in the function $\mathfrak{v}^{(l)}$), we assume the property that there exist $\mathbf{S}^{(l)}\in\mathbb R_\geq$ s.t. for all model parameters we have\\
\begin{flalign*}
&\Vert\mathfrak{f}^{(l)}(\uptheta^{(l)})-\mathfrak{v}^{(l)}(\check\uptheta^{(l)},\uptheta_\mathrm{m}^{(l)},\uptheta_\mathrm{v}^{(l)})\Vert^2 \leq (\mathbf{S}^{(l)})^2\Vert\uptheta^{(l)},\boldsymbol{0})^T-(\check\uptheta^{(l)},\uptheta_\mathrm{m}^{(l)},\uptheta_\mathrm{v}^{(l)})^T\Vert^2,\\
& \Vert(\uptheta^{(l)},\mathbf{0})^T-(\check\uptheta^{(l)},\uptheta_\mathrm{m}^{(l)},\uptheta_\mathrm{v})^T\Vert^2\neq 0,\quad  \forall (\uptheta^{(l)},\mathbf{0})^T,(\check\uptheta^{(l)},\uptheta_\mathrm{m}^{(l)},\uptheta_\mathrm{v}^{(l)})^T.
\end{flalign*}
where we expand the input-space of $\mathfrak{f}^{(l)}$ naturally with zeros, where we have variance parameters and where
\begin{flushleft}
$\uptheta^{(1)}_\mathrm{m}=\begin{pmatrix}\mathbf{m}^{(1)}\\\boldsymbol\upalpha^{(1)}\\\boldsymbol\upmu^{(1)}_{\mathsf{k}}\\\mathbf{\bar{x}}_{\mathsf{k}}\end{pmatrix}$, $\uptheta^{(l)}_\mathrm{m}=\begin{pmatrix}\mathbf{m}^{(l)}\\\boldsymbol\upalpha^{(l)}\\\boldsymbol\upmu^{(l)}_{\mathsf{k}}\\\boldsymbol\upmu^{(l-1)}_{\mathsf{k}}\end{pmatrix}$, $\uptheta^{(L+1)}_\mathrm{m}=\begin{pmatrix}\mathbf{m}^{(L+1)}\\\boldsymbol\upalpha^{(L+1)}\\\boldsymbol\upmu^{(L)}_{\mathsf{k}}\end{pmatrix}$,\\
 $\uptheta_\mathrm{v}^{(1)}=\begin{pmatrix}{\boldsymbol\upbeta}^{(1)}\\{\boldsymbol\uplambda}^{(1)}_{\mathsf{k}}\end{pmatrix}$, $\uptheta_\mathrm{v}^{(l)}=\begin{pmatrix}{\boldsymbol\upbeta}^{(l)}\\\boldsymbol\uplambda^{(l)}_{\mathsf{k}}\\\boldsymbol\uplambda^{(l-1)}_{\mathsf{k}}\end{pmatrix}$, $\uptheta_\mathrm{v}^{(L+1)}=\begin{pmatrix}{\boldsymbol\upbeta}^{(L+1)}\\\boldsymbol\uplambda^{(L)}_{\mathsf{k}}\end{pmatrix}$, for $l=2\dots,L$\\

$\check\uptheta^{(l)}_\mathfrak{v}=\begin{pmatrix}M\\(\sigma_{\text{power}}^{(l)})_\mathfrak{v}\\\frac{\text{diag}(\mathfrak{L}^{(l)}_\mathfrak{v})}{2\pi}\\\mathbf{p}^{(l)}_\mathfrak{v}\\\mathbf{u}^{(l)}_\mathfrak{v}\\\mathbf{b}^{(l)}_\mathfrak{v}\end{pmatrix}$, $\check\uptheta^{(l)}_\mathfrak{f}=\begin{pmatrix}M\\(\sigma_{\text{power}}^{(l)})_\mathfrak{f}\\\frac{\text{diag}(\mathfrak{L}^{(l)}_\mathfrak{f})}{2\pi}\\\mathbf{p}^{(l)}_\mathfrak{f}\\\mathbf{u}^{(l)}_\mathfrak{f}\\\mathbf{b}^{(l)}_\mathfrak{f} \end{pmatrix}$, $\hat\uptheta^{(l)}=\begin{pmatrix}\mathbf{a}^{(l)}\\\mathbf{z}^{(l)}\\\mathbf{x}^{(l)}_{\mathsf{k}}\end{pmatrix}$, for $l=1\dots,L+1$,
\end{flushleft}
In $\uptheta^{(1)}_\mathrm{m}$, $\hat\uptheta^{(1)}$ the $\mathbf{\bar{x}}_{\mathsf{k}}$ and in $\check\uptheta^{(l)}_\mathfrak{v}$, $\check\uptheta^{(l)}_\mathfrak{f}$ and $M$ argument cancels, so we are independent of this argument on the right hand side.\\
For the case when we choose $(\uptheta^{(l)}_\mathrm{m},\check\uptheta^{(l)}_\mathfrak{v})=(\hat\uptheta^{(l)},\check\uptheta^{(l)}_\mathfrak{f})$ only the variance parameters are naturally involved. If additionally $\uptheta_\mathrm{v}^{(1)}=\boldsymbol{0}$, $\mathfrak{v}^{(l)},\mathfrak{f}^{(l)}$ coincide.\\
We further assume the functions to be bounded with $\pm\frac{\mathbf{S}^{(l)}}{\updelta^{(l)}}$, which directly follows from the \textbf{(A1)} condition and additionally with $\updelta^{(l)}\in\mathbb R$ and
\begin{flalign*}
\Vert(\theta^{(l)},\boldsymbol{0})-(\check\theta^{(l)},\theta_\mathrm{m}^{(l)},\theta_\mathrm{v}^{(l)})\Vert^2 \leq (\updelta^{(l)})^2,\numberthis\label{kappa}
\end{flalign*}
by making an estimation up- and downwards.
\\
\noindent We go on showing the \textsl{quadratic-form-Gaussian} loss property for the case $\mathcal{L}^1(\lambda)$.\\
Therefore we assume $\underset{\substack{{\hat\theta \sim P_{\scalebox{.8}{\text{\tiny{REV}}}}}\\{\boldsymbol{y}\sim P_\mathsf{K}}}}{\mathbf{E}} = \underset{\substack{{\boldsymbol{y}\sim P_\mathsf{K}}\\{\hat\theta \sim P_{\scalebox{.8}{\text{\tiny{REV}}}}}}}{\mathbf{E}}$, Fubinis Theorem \textbf{(A2)}. We can also fill in assumption \textbf{(A4)} for the samples $\mathrm{y}'$ and $\mathrm{h}^{(l)}$, then we come to
\begin{align*}
\vartheta_{\mathcal{V}}(\lambda)&=\Psi^{\ell}(\lambda,1)\\
&=\log\left(\underset{\substack{{\boldsymbol{y}'\sim P_\mathsf{K}}\\{\hat\theta \sim P_{\scalebox{.8}{\text{\tiny{REV}}}}}}}{\mathbf{E}}\left[\exp\left(\lambda\left(\sum\nolimits_{l=1}^{L+1}\frac{(\mathbf{m}^{(l)}_*)^T(\Psi^{(l)}_{2,*}-(\Psi^{(l)}_{1,*})^T\Psi^{(l)}_{1,*})\mathbf{m}^{(l)}_*}{2(\sigma_{\text{noise}}^{(l)})^2}\right.\right.\right.\right.\\
&\left.\left.\left.\left.+\text{tr}\left(\frac{\Psi^{(l)}_{2,*}\mathbf{s}^{(l)}_*}{2(\sigma_{\text{noise}}^{(l)})^2}\right)+\frac{\mathsf{K}(\sigma_{\text{noise}*}^{(l)})^2}{2(\sigma_{\text{noise}}^{(l)})^2}+\frac{(\Psi^{(l)}_*\mathbf{m}^{(l)}_*-\Phi^{(l)}\mathbf{a}^{(l)})^T(\Psi^{(l)}_*\mathbf{m}^{(l)}_*-\Phi^{(l)}\mathbf{a}^{(l)})}{2(\sigma_{\text{noise}}^{(l)})^2}\right.\right.\right.\right.\\
&\left.\left.\left.\left.-\frac{\left(\Psi^{(L+1)}_*\mathbf{m}^{(L+1)}_*+\boldsymbol\epsilon^{\mathbf{y}}_*- \Phi^{(L+1)}\mathbf{a}^{(L+1)}\right)^T}{{2(\sigma_{\text{noise}}^{(L+1)})^2}}\frac{\left(\Psi^{(L+1)}_*\mathbf{m}^{(L+1)}_*+\boldsymbol\epsilon^{\mathbf{y}}_*- \Phi^{(L+1)}\mathbf{a}^{(L+1)}\right)}{{2(\sigma_{\text{noise}}^{(L+1)})^2}}\right.\right.\right.\right.\\
&\left.\left.\left.\left.-\sum\nolimits_{l=2}^{L} \frac{\left(\Psi^{(l)}_*\mathbf{m}^{(l)}_*+\boldsymbol\epsilon^{\mathbf{h}^{(l)}}_*- \Phi^{(l)}\mathbf{a}^{(l)}\right)^T}{2(\sigma_{\text{noise}}^{(l)})^2}\frac{\left(\Psi^{(l)}_*\mathbf{m}^{(l)}_*+\boldsymbol\epsilon^{\mathbf{h}^{(l)}}_*- \Phi^{(l)}\mathbf{a}^{(l)}\right)}{2(\sigma_{\text{noise}}^{(l)})^2}\right)\right)\right]\right)
\end{align*}
Calculating the binomials we arrive at
\begin{align*}
&=\log\left(\underset{\substack{{\boldsymbol{y}'\sim P_\mathsf{K}}\\{\hat\theta \sim P_{\scalebox{.8}{\text{\tiny{REV}}}}}}}{\mathbf{E}}\left[\exp\left(\lambda\left(\sum\nolimits_{l=1}^{L+1}\frac{(\mathbf{m}^{(l)}_*)^T(\Psi^{(l)}_{2,*}-(\Psi^{(l)}_{1,*})^T\Psi^{(l)}_{1,*})\mathbf{m}^{(l)}_*}{2(\sigma_{\text{noise}}^{(l)})^2}\right.\right.\right.\right.\\
&\left.\left.\left.\left.+\text{tr}\left(\frac{\Psi^{(l)}_{2,*}\mathbf{s}^{(l)}_*}{2(\sigma_{\text{noise}}^{(l)})^2}\right)+\frac{\mathsf{K}(\sigma_{\text{noise}*}^{(l)})^2}{2(\sigma_{\text{noise}}^{(l)})^2}\right.\right.\right.\right.\\
&\left.\left.\left.\left.-\frac{(\boldsymbol\epsilon^{\mathbf{y}}_*)^2+2\boldsymbol\epsilon^{\mathbf{y}}_*\left( \Psi^{(L+1)}_*\mathbf{m}^{(L+1)}_*- \Phi^{(L+1)}\mathbf{a}^{(L+1)}\right)}{2(\sigma_{\text{noise}}^{(L+1)})^2} \right.\right.\right.\right.\\
&\left.\left.\left.\left.-\sum\nolimits_{l=2}^{L+1}\frac{(\boldsymbol\epsilon^{\mathbf{h}^{(l)}}_*)^2+2\boldsymbol\epsilon^{\mathbf{h}^{(l)}}_*\left(\Psi^{(l)}_*\mathbf{m}^{(l)}_*- \Phi^{(l)}\mathbf{a}^{(l)}\right)}{2(\sigma_{\text{noise}}^{(l)})^2}\right)\right)\right]\right).
	\end{align*}
\\We write for DRGP-(V)SS, DRGP-Nyström
\begin{align*}
{\displaystyle 1\!\!1_\mathsf{K}^T}\textbf{Var}\left[\boldsymbol{h}^{(l)}\right] = (\mathbf{m}^{(l)}_*)^T(\Psi^{(l)}_{2,*}-(\Psi^{(l)}_{1,*})^T\Psi^{(l)}_{1,*})\mathbf{m}^{(l)}_*+\text{tr}\left(\Psi^{(l)}_{2,*}\mathbf{s}^{(l)}_*\right)+\mathsf{K}(\sigma_{\text{noise}*}^{(l)})^2,
\end{align*}
and for DRGP-SS, DRGP-Nyström we have
\begin{align*}
\textbf{Cov}\left[\boldsymbol{h}^{(l)}\right] = \Psi^{(l)}_{1,*}\mathbf{s}^{(l)}_*(\Psi^{(l)}_{1,*})^T\in\mathbb R^{\mathsf{K}\times \mathsf{K}}, \mathsf{k}\neq \hat{\mathsf{k}}, \mathsf{k},\hat{\mathsf{k}}=1,\dots,\mathsf{K},
\end{align*}
and for DRGP-VSS we have
\begin{align*}
\textbf{Cov}\left[{h}^{(l)}_k,{h}^{(l)}_{k'}\right] = (\mathbf{m}^{(l)}_*)^T(\Psi^{(l)}_{\mathsf{k},\hat{\mathsf{k}},*}-(\psi^{(l)}_{1,*,\mathsf{k}})^T\psi^{(l)}_{1,*,\hat{\mathsf{k}}})\mathbf{m}^{(l)}_*+\text{tr}\left(\Psi^{(l)}_{\mathsf{k},\hat{\mathsf{k}},*}\mathbf{s}^{(l)}_*\right), \mathsf{k}\neq \hat{\mathsf{k}}, \hat{\mathsf{k}}=1,\dots,\mathsf{K},
\end{align*}
where $\Psi^{(l)}_{\mathsf{k},\hat{\mathsf{k}},*}=\sum_{m=1}^M(\psi)^{m,(l)}_{\mathsf{k},\hat{\mathsf{k}},*}$, which is $(\psi^{(l)}_{1,*,\mathsf{k}})^T\psi^{(l)}_{1,*,\hat{\mathsf{k}}}-\text{diag}[(\psi^{(l)}_{1,*,\mathsf{k}})^T]+ D$ where $D$ is a diagonal matrix given by:
\begin{align*}
D&=\text{diag}\left[\frac{(2\pi)^Q\upsigma_{\text{power}}^2\prod\limits_{q=1}^Q\left(\frac{\mathrm{l}_q^2}{\boldsymbol\uplambda_{n_q}+\boldsymbol\uplambda_{n'_q}}\right)}{M}(e^{-\frac{1}{2}(\bar{\boldsymbol\upmu}_{nm}-\bar{\boldsymbol\upmu}_{n'm})^T ((2\pi)^{-2}\mathfrak{L}^2(\boldsymbol\uplambda_n+\boldsymbol\uplambda_{n'})^{-1})(\bar{\boldsymbol\upmu}_{nm}-\bar{\boldsymbol\upmu}_{n'm})}\right.\\
&\left.\cos(-\mathbf{p}_m^T(\bar{\boldsymbol\upmu}_{nm}-\bar{\boldsymbol\upmu}_{n'm}))+e^{-\frac{1}{2}(\bar{\boldsymbol\upmu}_{nm}+\bar{\boldsymbol\upmu}_{n'm})^T ((2\pi)^{-2}\mathfrak{L}^2(\boldsymbol\uplambda_n+\boldsymbol\uplambda_{n'})^{-1})(\bar{\boldsymbol\upmu}_{nm}+\bar{\boldsymbol\upmu}_{n'm})}\right.\\
&\left.\cos(-\mathbf{p}_m^T(\bar{\boldsymbol\upmu}_{nm}+\bar{\boldsymbol\upmu}_{n'm})))\right]_{m=1}^M
\end{align*}
Calculating now the moment-generating functions for $\boldsymbol\epsilon^{\mathbf{y}}_*$, $\boldsymbol\epsilon^{\mathbf{h}^{(l)}}_*$, see~\cite[Section 3.2]{olkin1992quadratic}, 
we come to 
\begin{align*}
&=\sum\nolimits_{l=1}^{L+1}\frac{\lambda {\displaystyle 1\!\!1_\mathsf{K}^T}\textbf{Var}\left[\boldsymbol{h}^{(l)}\right]}{2(\sigma_{\text{noise}}^{(l)})^2}-\frac{1}{2}\log\left(|I_\mathsf{K}+2\frac{\lambda\textbf{Cov}\left[\boldsymbol{h}^{(l)}\right]}{2(\sigma_{\text{noise}}^{(l)})^2}|\right)\\
&+\log\left(\underset{{\theta \sim P_{\scalebox{.8}{\text{\tiny{REV}}}}}}{\mathbf{E}}\left[\exp\left(\frac{1}{2}\left(\frac{\left(\Psi^{(l)}_*\mathbf{m}^{(l)}_*- \Phi^{(l)}\mathbf{a}^{(l)}\right)^T\lambda}{(\sigma_{\text{noise}}^{(l)})^2} \textbf{Cov}\left[\boldsymbol{h}^{(l)}\right]\right)^T\right.\right.\right.\\
&\left.\left.\left.\left(I_\mathsf{K}+\frac{\lambda \textbf{Cov}\left[\boldsymbol{h}^{(l)}\right]}{(\sigma_{\text{noise}}^{(l)})^2}\right)^{-1}\frac{\lambda\left(\Psi^{(l)}_*\mathbf{m}^{(l)}_*- \Phi^{(l)}\mathbf{a}^{(l)}\right)}{(\sigma_{\text{noise}}^{(l)})^2}\right)\right]\right),
\end{align*}
and moreover, we can rewrite
\begin{align*}
\textbf{Cov}\left[\boldsymbol{h}^{(l)}\right]\left(I_\mathsf{K}+\frac{\lambda\textbf{Cov}\left[\boldsymbol{h}^{(l)}\right]}{(\sigma_{\text{noise}}^{(l)})^2}\right)^{-1}=\left(\textbf{Cov}\left[\boldsymbol{h}^{(l)}\right]^{-1}+\frac{\lambda I_\mathsf{K}}{(\sigma_{\text{noise}}^{(l)})^2}\right)^{-1}.\numberthis\label{L10}
\end{align*}
and further, as the expression in the exponential is always positive, because of positive definiteness of the matrix in Equation~(\ref{L10}), $e^{(.)}$ is always greater than 1 and therefore $\log(.)$ is always non-negative. 
We further see, that the matrix in Equation~(\ref{L10}) is also real, symmetric and therefore defines a scalar product. We therefore can apply the stochastic Lipschitz condition \textbf{(A1)} followed by the condition on the input-space (overall a bound condition on the mean functions)
\begin{align*}
&\leq\sum\nolimits_{l=1}^{L+1}\frac{\lambda {\displaystyle 1\!\!1_\mathsf{K}^T}\textbf{Var}\left[\boldsymbol{h}^{(l)}\right]}{2(\sigma_{\text{noise}}^{(l)})^2}-\sum\nolimits_{\mathsf{k}=1}^\mathsf{K}\frac{1}{2}\log\left(|I_\mathsf{K}+\frac{\lambda\textbf{Cov}\left[\boldsymbol{h}^{(l)}\right]}{(\sigma_{\text{noise}}^{(l)})^2}|\right)\\
&+\frac{1}{2}\left(\left[\frac{\mathbf{S}^{(l)}}{\updelta}\frac{\lambda}{(\sigma_{\text{noise}}^{(l)})^2} \right]_{\mathsf{k}=1}^\mathsf{K}\right)^T\left(\textbf{Cov}\left[\boldsymbol{h}^{(l)}\right]\right)\left(I_\mathsf{K}+\frac{\lambda \textbf{Cov}\left[\boldsymbol{h}^{(l)}\right]}{(\sigma_{\text{noise}}^{(l)})^2}\right)^{-1}\left[\frac{\mathbf{S}^{(l)}}{\updelta}\frac{\lambda}{(\sigma_{\text{noise}}^{(l)})^2} \right]_{\mathsf{k}=1}^\mathsf{K},
\end{align*}
which we denote with $\mathbb{L}(\lambda)$ and which proves the property of the loss function. Further one can show, that this expression is always positive.\\
\\
Assuming now $i=1,\dots,N$ iid observations of the random variable $\mathcal{Q}(\mathbfcal{I})$, denoted with $l_i$, defining $\mathcal{E}^{\mathbf{y}^{N}}_*=\begin{pmatrix}\boldsymbol\epsilon^{\mathbf{y}^1}_*\dots\boldsymbol\epsilon^{\mathbf{y}^N}_*\end{pmatrix}^T$, $\mathcal{E}^{(\mathbf{h}^{(l)})^{N}}_*=\begin{pmatrix}\boldsymbol\epsilon^{\mathbf{h}^{(l),1}}_*\dots\boldsymbol\epsilon^{\mathbf{h}^{(l),N}}_*\end{pmatrix}^T$, then the expression for $\mathbb{L}(\lambda)$ becomes
\begin{align*}
\Psi^{\ell}(\lambda,N)&= \log\left(\underset{\substack{{\mathbbmsl{D}'\sim (P_\mathsf{K})^N}\\{\hat\theta \sim P_{\scalebox{.8}{\text{\tiny{REV}}}}}}}{\mathbf{E}}\left[e^{\frac{\lambda}{N}\sum\nolimits_{i=1}^N l_i}\right]\right)\\
&=\dots\\
&=\log\left(\underset{\substack{{\mathbbmsl{D}'\sim (P_\mathsf{K})^N}\\{\hat\theta \sim P_{\scalebox{.8}{\text{\tiny{REV}}}}}}}{\mathbf{E}}\left[\exp\left(\frac{\lambda}{N}\left(N\sum\nolimits_{l=1}^{L+1} \frac{(\mathbf{m}^{(l)}_*)^T(\Psi^{(l)}_{2,*}-(\Psi^{(l)}_{1,*})^T\Psi^{(l)}_{1,*})\mathbf{m}^{(l)}_*}{2(\sigma_{\text{noise}}^{(l)})^2}\right.\right.\right.\right.\\
&\left.\left.\left.\left.+\text{tr}\left(\frac{\Psi^{(l)}_{2,*}\mathbf{s}^{(l)}_*}{2(\sigma_{\text{noise}}^{(l)})^2}\right)+\frac{\mathsf{K}(\sigma_{\text{noise}*}^{(l)})^2 }{2(\sigma_{\text{noise}}^{(l)})^2}\right.\right.\right.\right.\\
&\left.\left.\left.-\frac{\text{tr}(\mathcal{E}^{(\mathbf{y})^N}_*(\mathcal{E}^{(\mathbf{y})^N}_*)^T)+2\left(\Psi^{(L+1)}_*\mathbf{m}^{(L+1)}_*- \Phi^{(L+1)}\mathbf{a}^{(L+1)}\right)(\mathcal{E}^{(\mathbf{y})^N}_*)^T{\displaystyle 1\!\!1_{N}}}{2(\sigma_{\text{noise}}^{(L+1)})^2}\right.\right.\right.\\
&\left.\left.\left.-\sum\nolimits_{l=1}^{L}\frac{\text{tr}(\mathcal{E}^{(\mathbf{h}^{(l)})^N}_*(\mathcal{E}^{(\mathbf{h}^{(l)})^N}_*)^T)+2\left(\Psi^{(l)}_*\mathbf{m}^{(l)}_*- \Phi^{(l)}\mathbf{a}^{(l)}\right)(\mathcal{E}^{(\mathbf{h}^{(l)})^N}_*)^T{\displaystyle 1\!\!1_{N}}}{2(\sigma_{\text{noise}}^{(l)})^2}
\right)\right]\right)
\end{align*}
Because of the independence assumption, this can now be calculated as before to
\begin{align*}
&\leq \sum\nolimits_{l=1}^{L+1}\frac{\lambda {\displaystyle 1\!\!1_\mathsf{K}^T}\textbf{Var}\left[\boldsymbol{h}^{(l)}\right]}{2(\sigma_{\text{noise}}^{(l)})^2}-N\frac{1}{2}\log\left(|I_\mathsf{K}+\frac{1}{N}\frac{\lambda \textbf{Cov}\left[\boldsymbol{h}^{(l)}\right]}{(\sigma_{\text{noise}}^{(l)})^2}|\right)\numberthis\label{L2}\\
&+N\frac{1}{2}\left(\left[\frac{1}{N}\frac{\mathbf{S}^{(l)}}{\updelta}\frac{\lambda }{(\sigma_{\text{noise}}^{(l)})^2}\right]_{\mathsf{k}=1}^\mathsf{K}\right)^T\left(\textbf{Cov}\left[\boldsymbol{h}^{(l)}\right]\right)\\
&\left(I_\mathsf{K}+\frac{1}{N}\frac{\lambda \textbf{Cov}\left[\boldsymbol{h}^{(l)}\right]}{(\sigma_{\text{noise}}^{(l)})^2}\right)^{-1}\left[\frac{1}{N}\frac{\mathbf{S}^{(l)}}{\updelta}\frac{\lambda }{(\sigma_{\text{noise}}^{(l)})^2}\right]_{\mathsf{k}=1}^\mathsf{K},
\end{align*}
which we denote with $\mathcal{L}(\lambda)$.\\
Then with either $\lambda=N$ or $\lambda=\sqrt{N}$ we come to
\begin{align*}
\Psi^{\ell}(\lambda,N)&= \log\left(\underset{\substack{{\mathbbmsl{D}'\sim P_\mathsf{K}}\\{\hat\theta \sim P_{\scalebox{.8}{\text{\tiny{REV}}}}}}}{\mathbf{E}}\left[e^{\frac{\lambda}{N}\sum\nolimits_{i=1}^N l_i}\right]\right)\\
&\leq\dots\\
&\stackrel{\lambda=N}=\mathcal{L}\left(N\right)\\
&\stackrel{\lambda=\sqrt{N}}=\mathcal{L}\left(\sqrt{N}\right)
\end{align*}
\\
This is different to the derivation of ~\cite{germain2016pac} (see also Annotation $^3$ in~\cite{NIPS20177100}). We have
\begin{align*}
\frac{\mathcal{L}(\sqrt{N})}{\sqrt{N}}&=\sum\nolimits_{l=1}^{L+1}\frac{{\displaystyle 1\!\!1_\mathsf{K}^T}\textbf{Var}\left[\boldsymbol{h}^{(l)}\right]}{2(\sigma_{\text{noise}}^{(l)})^2}-\frac{\sqrt{N}}{2}\log\left(|I_\mathsf{K}+\frac{
 \textbf{Cov}\left[\boldsymbol{h}^{(l)}\right]}{\sqrt{N}(\sigma_{\text{noise}}^{(l)})^2}|\right)\\
&+\frac{1}{2\sqrt{N}}\left(\left[\frac{\mathbf{S}^{(l)}}{\updelta}\frac{1}{(\sigma_{\text{noise}}^{(l)})^2}\right]_{\mathsf{k}=1}^\mathsf{K}\right)^T\left(\textbf{Cov}\left[\boldsymbol{h}^{(l)}\right]\right)\\
&\left(I_\mathsf{K}+\frac{\textbf{Cov}\left[\boldsymbol{h}^{(l)}\right]}{\sqrt{N}(\sigma_{\text{noise}}^{(l)})^2}\right)^{-1}\left[\frac{\mathbf{S}^{(l)}}{\updelta}\frac{1}{(\sigma_{\text{noise}}^{(l)})^2}\right]_{\mathsf{k}=1}^\mathsf{K},
\end{align*}
\begin{align*}
\frac{\mathcal{L}(N)}{N}&=\sum\nolimits_{l=1}^{L+1}\frac{{\displaystyle 1\!\!1_\mathsf{K}^T}\textbf{Var}\left[\boldsymbol{h}^{(l)}\right]}{2(\sigma_{\text{noise}}^{(l)})^2}-\frac{1}{2}\log\left(|I_\mathsf{K}+\frac{
 \textbf{Cov}\left[\boldsymbol{h}^{(l)}\right]}{(\sigma_{\text{noise}}^{(l)})^2}|\right)\\
&+\frac{1}{2}\left(\left[\frac{\mathbf{S}^{(l)}}{\updelta}\frac{1}{(\sigma_{\text{noise}}^{(l)})^2}\right]_{\mathsf{k}=1}^\mathsf{K}\right)^T\left(\textbf{Cov}\left[\boldsymbol{h}^{(l)}\right]\right)\\
&\left(I_\mathsf{K}+\frac{\textbf{Cov}\left[\boldsymbol{h}^{(l)}\right]}{(\sigma_{\text{noise}}^{(l)})^2}\right)^{-1}\left[\frac{\mathbf{S}^{(l)}}{\updelta}\frac{1}{(\sigma_{\text{noise}}^{(l)})^2}\right]_{\mathsf{k}=1}^\mathsf{K}.
\end{align*}
Convergence is achieved in the case $\lambda=\sqrt{N}$. This can be seen, because we assumed bounded covariance entries in \textbf{(A3)}. Therefore, Equation~(\ref{L10}) is well defined and exists. Therefore the last to rows convergence to zero for $N\to\infty$. Regarding the first row, as our models always reproduce positive definite, real, symmetric matrices, ${\displaystyle 1\!\!1_\mathsf{K}^T}\textbf{Var}\left[\boldsymbol{h}^{(l)}\right]=\text{tr}(\textbf{Cov}\left[\boldsymbol{h}^{(l)}\right])$ and $\text{sum}(\text{eig}(\textbf{Cov}\left[\boldsymbol{h}^{(l)}\right]))$ coincides. As $\sqrt{N}\log(1+\frac{x}{\sqrt{N}})$ converges to $x$ for $N\to\infty$, also the first row vanishes.

\newpage
\section{Proof of Theorem 5.}
\label{pacbayesian2}
{\bf Proof}\\
\\
The inequality in Theorem 4. can be rewritten with $\mathcal{L}_{ P_\mathsf{K}}^{\ell_{\text{nll}}}(\mathfrak{f}_{\bar{\uptheta}})=\text{arg}\min\nolimits_{\hat\uptheta\in\hat{\uptheta}}\left(\mathcal{L}_{P_\mathsf{K}}^{\ell_{\text{nll}}}(\mathfrak{f}_{\uptheta})\right)$ (Bayes risk constrained to our raw model with parameters restricted to the variational ones) to 
\begin{align*}
&\underset{{\hat\theta\sim Q_{\lambda}}}{\mathbf{E}}\left[\mathcal{L}_{P_\mathsf{K}}^{\ell_{\text{nll}}}(\mathfrak{f}_{\uptheta})\right]-\mathcal{L}_{P_\mathsf{K}}^{\ell_{\text{nll}}}(\mathfrak{f}_{\bar{\uptheta}})\leq\inf\limits_{Q_{\text{\scalebox{.8}{\tiny{REV}}}}\text{ in }\hat{\mathcal G}_{\boldsymbol\Theta}}\left(\underset{\hat\theta\sim Q_{\text{\scalebox{.8}{\tiny{REV}}}}}{\mathbf{E}}\left[\mathcal{L}_{P_\mathsf{K}}^{\ell_{\text{nll}}}(\mathfrak{f}_{\uptheta})-\mathcal{L}_{P_\mathsf{K}}^{\ell_{\text{nll}}}(\mathfrak{f}_{\bar{\uptheta}})\right]\right.\\
& \left.+\frac{2}{\lambda}\left(\mathbf{KL}(Q_{\text{\scalebox{.8}{\tiny{REV}}}}||P_{\text{\scalebox{.8}{\tiny{REV}}}})+\log\left(\frac{1}{\uptau}\right)+\Psi^{\ell}(\lambda,N)\right)\right).
\end{align*}
We show a bound on $\underset{\hat\theta\sim Q_{\text{\scalebox{.8}{\tiny{REV}}}}}{\mathbf{E}}\left[\mathcal{L}_{P_\mathsf{K}}^{\ell_{\text{nll}}}(\mathfrak{f}_{\bar{\uptheta}})-\mathcal{L}_{P_\mathsf{K}}^{\ell_{\text{nll}}}(\mathfrak{f}_{\uptheta})\right]$ with assumption \textbf{(A4)} similar to the procedure in Theorem 2. and 3., where we assume $\underset{\substack{{\hat\theta \sim Q_{\scalebox{.8}{\text{\tiny{REV}}}}}\\{\boldsymbol{y}\sim P_\mathsf{K}}}}{\mathbf{E}} = \underset{\substack{{\boldsymbol{y}\sim P_\mathsf{K}}\\{\hat\theta \sim Q_{\scalebox{.8}{\text{\tiny{REV}}}}}}}{\mathbf{E}}$, Fubinis Theorem \textbf{(A2)}. Then we come to 
\begin{align*}
&\underset{\hat\theta\sim Q_{\text{\scalebox{.8}{\tiny{REV}}}}}{\mathbf{E}}\left[\mathcal{L}_{P_{K}}^{\ell_{\text{nll}}}(\mathfrak{f}_{\uptheta})-\mathcal{L}_{P_{K}}^{\ell_{\text{nll}}}(\mathfrak{f}_{\bar{\uptheta}})\right]\\
&=\underset{\substack{{\boldsymbol{y}\sim P_\mathsf{K}}\\{\hat\theta \sim Q_{\scalebox{.8}{\text{\tiny{REV}}}}}}}{\mathbf{E}}\left[\ell_{\text{nll}}(\mathfrak{f}_{\uptheta},\mathbf{y})-\ell_{\text{nll}}(\mathfrak{f}_{\bar{\uptheta}},\mathbf{y})\right]\\
&=\underset{\substack{{\boldsymbol{y}\sim P_\mathsf{K}}\\{\hat\theta \sim Q_{\scalebox{.8}{\text{\tiny{REV}}}}}}}{\mathbf{E}}\left[-\left(\log(p(\mathbf{y}|\uptheta^{(L+1)},\mathbf{X}^{(L+1)}))+\sum\nolimits_{l=1}^{L}\log(p(\mathbf{h}^{(l)}|\uptheta^{(l)},\mathbf{X}^{(l)}))\right)\right.\\
&\left.+\log(p(\mathbf{y}|\uptheta^{(L+1)},\mathbf{\bar{X}}^{(L+1)}))+\sum\nolimits_{l=1}^{L}\log(p(\mathbf{\bar{h}}^{(l)}|\uptheta^{(l)},\mathbf{\bar{X}}^{(l)}))\right]\\
&=\underset{\substack{{\hat\theta \sim Q_{\scalebox{.8}{\text{\tiny{REV}}}}}\\{\boldsymbol{y}\sim P_\mathsf{K}}}}{\mathbf{E}}\left[\frac{(\mathbf{y}- \Phi^{(L+1)}\mathbf{a}^{(L+1)})^T(\mathbf{y}- \Phi^{(L+1)}\mathbf{a}^{(L+1)})}{2(\sigma_{\text{noise}}^{(L+1)})^2}+\mathsf{K}\frac{\log\left(2\pi(\sigma_{\text{noise}}^{(L+1)})^2\right)}{2}\right.\\
&\left.+\sum\nolimits_{l=1}^L \frac{(\mathbf{h}^{(l)}- \Phi^{(l)}\mathbf{a}^{(l)})^T(\mathbf{h}^{(l)}- \Phi^{(l)}\mathbf{a}^{(l)})}{2(\sigma_{\text{noise}}^{(l)})^2}+\mathsf{K}\frac{\log\left(2\pi(\sigma_{\text{noise}}^{(l)})^2\right)}{2}\right.\\
&\left.-\frac{(\mathbf{y}- \bar{\Phi}^{(L+1)}\mathbf{\bar{a}}^{(L+1)})^T(\mathbf{y}- \bar{\Phi}^{(L+1)}\mathbf{\bar{a}}^{(L+1)})}{2(\sigma_{\text{noise}}^{(L+1)})^2}-\mathsf{K}\frac{\log\left(2\pi(\sigma_{\text{noise}}^{(L+1)})^2\right)}{2}\right.\\
&\left.-\sum\nolimits_{l=1}^L \frac{(\mathbf{\bar{h}}^{(l)}- \bar{\Phi}^{(l)}\mathbf{\bar{a}}^{(l)})^T(\mathbf{\bar{h}}^{(l)}-\bar{\Phi}^{(l)}\mathbf{\bar{a}}^{(l)})}{2(\sigma_{\text{noise}}^{(l)})^2}+\mathsf{K}\frac{\log\left(2\pi(\sigma_{\text{noise}}^{(l)})^2\right)}{2}\right]
\end{align*}
\begin{align*}
&=\underset{\boldsymbol{y}\sim P_\mathsf{K}}{\mathbf{E}}\left[\frac{-2\mathbf{y}^T\Psi_1^{\text{(L+1)}}\mathbf{m}^{(L+1)}}{2(\upsigma_{\text{noise}}^{\text{(L+1)}})^2}+\frac{\text{tr}\left(\Psi_2^{\text{(L+1)}}(\mathbf{s}^{(L+1)} + \mathbf{m}^{(L+1)}(\mathbf{m}^{(L+1)})^T)\right)}{2(\upsigma_{\text{noise}}^{\text{(L+1)}})^2}\right.\\
&\left.+\sum\nolimits_{l=1}^{L}\frac{{\displaystyle 1\!\!1_\mathsf{K}^T}\boldsymbol\uplambda^{(l)}+ (\boldsymbol\upmu^{(l)})^T\boldsymbol\upmu^{(l)}-2\upmu^{(l)}\Psi_1^{(l)}\mathbf{m}^{(l)}}{2(\upsigma_{\text{noise}}^{(l)})^2}+\frac{\text{tr}\left(\Psi_2^{(l)}(\mathbf{s}^{(l)} + \mathbf{m}^{(l)}(\mathbf{m}^{(l)})^T)\right)}{2(\upsigma_{\text{noise}}^{(l)})^2}\right.\\
&\left.- \frac{-2\mathbf{y}^T\bar{\Phi}^{(L+1)}\mathbf{\bar{a}}^{(L+1)}+(\mathbf{\bar{a}}^{(L+1)})^T(\bar{\Phi}^{(L+1)})^T\bar{\Phi}^{(L+1)}\mathbf{\bar{a}}^{(L+1)})}{2(\sigma_{\text{noise}}^{(L+1)})^2}\right.\\
&\left.-\sum\nolimits_{l=1}^L \frac{(\mathbf{\bar{h}}^{(l)})^T\mathbf{\bar{h}}^{(l)}-2(\mathbf{\bar{h}}^{(l)})^T\bar{\Phi}^{(l)}\mathbf{\bar{a}}^{(l)})^T+(\mathbf{\bar{a}}^{(l)})^T(\bar{\Phi}^{(l)})^T\bar{\Phi}^{(l)}\mathbf{\bar{a}}^{(l)})}{2(\sigma_{\text{noise}}^{(l)})^2}\right]
\end{align*}
After setting $\mathbf{\bar{h}}^{(l)}=\boldsymbol\upmu^{(l)}$ and calculating the expectation we come to
\begin{align*}
&=\frac{-2(\mathbf{m}^{(l)}_*)^T(\Psi_{1,*}^{(l)})^T\Psi_1^{\text{(L+1)}}\mathbf{m}^{(L+1)}}{2(\upsigma_{\text{noise}}^{\text{(L+1)}})^2}+\frac{\text{tr}\left(\Psi_2^{\text{(L+1)}}(\mathbf{s}^{(L+1)} + \mathbf{m}^{(L+1)}(\mathbf{m}^{(L+1)})^T)\right)}{2(\upsigma_{\text{noise}}^{\text{(L+1)}})^2}\\
&+\sum\nolimits_{l=1}^{L}\frac{{\displaystyle 1\!\!1_\mathsf{K}^T}\boldsymbol\uplambda^{(l)}-2\boldsymbol\upmu^{(l)}\Psi_1^{(l)}\mathbf{m}^{(l)}}{2(\upsigma_{\text{noise}}^{(l)})^2}+\frac{\text{tr}\left(\Psi_2^{(l)}(\mathbf{s}^{(l)} + \mathbf{m}^{(l)}(\mathbf{m}^{(l)})^T)\right)}{2(\upsigma_{\text{noise}}^{(l)})^2}\\
&- \frac{-2(\mathbf{m}^{(l)}_*)^T(\Psi_{1,*}^{(l)})^T\hat{\Phi}^{(L+1)}\mathbf{\bar{a}}^{(L+1)})^T+(\mathbf{\bar{a}}^{(L+1)})^T(\hat{\Phi}^{(L+1)})^T\hat{\Phi}^{(L+1)}\mathbf{\bar{a}}^{(L+1)})}{2(\sigma_{\text{noise}}^{(L+1)})^2}\\
&-\sum\nolimits_{l=1}^L \frac{-2(\boldsymbol\upmu^{(l)})^T\hat{\Phi}^{(l)}\mathbf{\bar{a}}^{(l)})^T+(\mathbf{\bar{a}}^{(l)})^T(\hat{\Phi}^{(l)})^T\hat{\Phi}^{(l)}\mathbf{\bar{a}}^{(l)})}{2(\sigma_{\text{noise}}^{(l)})^2},
\end{align*}
where we can replace $\mathbf{m}^{(l)}$ (respective *) throughout these calculations with $(A^{(l)})^{-1}\left(\hat{\Psi}_1^{(l)}\right)^T\upmu^{(l)}$ for $l=1,\dots,L$ (respective $(A^{(L+1)})^{-1}\left(\hat{\Psi}_1^{(L+1)}\right)^T\mathbf{Y}$) and $\mathbf{s}^{(l)}$ with $(\sigma_{\text{noise}}^{(l)})^2(A^{(l)})^{-1}$ for $l=1,\dots,L+1$. $\hat{\Phi}$ is again the matrix filled in with the optimal parameters where no prior assumption exists.\\
Replacing $\upmu^{(l)}$ with $\Psi_1^{(l)}\mathbf{m}^{(l)}$ and making further simplifications
\begin{align*}
&=\frac{-2(\mathbf{m}^{(l)}_*)^T(\Psi_{1,*}^{(l)})^T\Psi_1^{\text{(L+1)}}\mathbf{m}^{(L+1)}}{2(\upsigma_{\text{noise}}^{\text{(L+1)}})^2}+\frac{\text{tr}\left(\Psi_2^{\text{(L+1)}}(\mathbf{s}^{(L+1)} + \mathbf{m}^{(L+1)}(\mathbf{m}^{(L+1)})^T)\right)}{2(\upsigma_{\text{noise}}^{\text{(L+1)}})^2}\\
&- \frac{-2(\mathbf{m}^{(l)}_*)^T(\Psi_{1,*}^{(l)})^T\bar{\Phi}^{(L+1)}\mathbf{\bar{a}}^{(L+1)})^T+(\mathbf{\bar{a}}^{(L+1)})^T(\bar{\Phi}^{(L+1)})^T\bar{\Phi}^{(L+1)}\mathbf{\bar{a}}^{(L+1)})}{2(\sigma_{\text{noise}}^{(L+1)})^2}\\
&+\sum\nolimits_{l=1}^{L}\frac{{\displaystyle 1\!\!1_\mathsf{K}^T}\boldsymbol\uplambda^{(l)}}{2(\upsigma_{\text{noise}}^{(l)})^2}+\frac{{\displaystyle 1\!\!1_\mathsf{K}^T}\textbf{Var}\left[\boldsymbol{f}^{(l)}_{\mathbf{X}^{(l)},\uptheta}\right]}{2(\sigma_{\text{noise}}^{(l)})^2}-\frac{(\Psi_{1}^{(l)}\mathbf{m}^{(l)}-\hat{\Phi}^{(l)}\mathbf{\bar{a}}^{(l)})^T(\Psi_{1}^{(l)}\mathbf{m}^{(l)}-\hat{\Phi}^{(l)}\mathbf{\bar{a}}^{(l)})}{2(\sigma_{\text{noise}*}^{(l)})^2}.
\end{align*}
\\
\textbf{(D1)} from the beginning also holds here.\\
\\
Setting now ones again the variational parameters such that the mean is equal to the optimal values $\mathbf{m}^{(l)}_*,\boldsymbol\alpha^{(l)}_*,\boldsymbol\upmu^{(l)}_*,\mathbf{s}^{(l)}_*,\boldsymbol\upbeta^{(l)}_*,\boldsymbol\uplambda^{(l)}_*$, $\sigma_{\text{noise}*}^{(l)}$, this expression becomes insightful in terms of the expressions which are involved
\begin{align*}
&=\sum\nolimits_{l=1}^{L}\frac{{\displaystyle 1\!\!1_\mathsf{K}^T}\boldsymbol\uplambda_*^{(l)}}{2(\upsigma_{\text{noise}*}^{(l)})^2}+\sum\nolimits_{l=1}^{L+1}\frac{{\displaystyle 1\!\!1_\mathsf{K}^T}\textbf{Var}\left[\boldsymbol{f}^{(l)}_{\mathbf{X}^{(l)},\uptheta_*}\right]}{2(\sigma_{\text{noise}*}^{(l)})^2}-\frac{(\Psi_{1,*}^{(l)}\mathbf{m}^{(l)}_*-\hat{\Phi}^{(l)}\mathbf{\bar{a}}^{(l)})^T(\Psi_{1,*}^{(l)}\mathbf{m}^{(l)}_{*}-\hat{\Phi}^{(l)}\mathbf{\bar{a}}^{(l)})}{2(\sigma_{\text{noise}*}^{(l)})^2}.
\end{align*}
Setting now $\lambda=\sqrt{N}$ and deleting the last expression, as it is always positive, we derive at
\begin{align*}
&\underset{{\hat\theta\sim Q_{\lambda}}}{\mathbf{E}}\left[\mathcal{L}_{P_\mathsf{K}}^{\ell_{\text{nll}}}(\mathbf{f}_{\uptheta})\right]-\mathcal{L}_{P_\mathsf{K}}^{\ell_{\text{nll}}}(\mathbf{f}_{\uptheta_*})\leq \sum\nolimits_{l=1}^{L}\frac{{\displaystyle 1\!\!1_\mathsf{K}^T}\boldsymbol\uplambda_*^{(l)}}{2(\upsigma_{\text{noise*}}^{(l)})^2}+ \sum\nolimits_{l=1}^{L+1}\frac{{\displaystyle 1\!\!1_\mathsf{K}^T}\textbf{Var}\left[\boldsymbol{f}^{(l)}_{\mathbf{X}^{(l)},\uptheta_*}\right]}{2(\sigma_{\text{noise}*}^{(l)})^2}\\
&+\frac{2}{\sqrt{N}}\left(\mathbf{KL}(Q_{*}||P_{\text{\scalebox{.8}{\tiny{REV}}}})+\log\left(\frac{1}{\uptau}\right)+\mathcal{L}(\sqrt{N})\right).
\end{align*}
For $\lambda=\sqrt{N}$ we have convergence to $\sum\nolimits_{l=1}^{L+1}\frac{{\displaystyle 1\!\!1_\mathsf{K}^T}\left(\boldsymbol\uplambda_*^{(l)}+\textbf{Var}\left[\boldsymbol{f}^{(l)}_{\mathbf{X}^{(l)},\uptheta_*}\right]\right)}{2(\sigma_{\text{noise}*}^{(l)})^2}$ with $\boldsymbol\uplambda_*^{(L+1)}=0$.
\newpage
\section{Some minor proofs of the extensions}
\label{extensions}
We have for the extension case\\
\\
$R=\frac{\mathbf{KL}(Q_{\text{\scalebox{.8}{\tiny{REV}}}}||P_{\text{\scalebox{.8}{\tiny{REV}}}})+\log\left(\frac{1}{\uptau}\right)+\log\left(\underset{\substack{{\hat\theta \sim P_{\scalebox{.8}{\text{\tiny{PAC}}}}},{\mathbbmsl{D}'\sim (P_\mathsf{K})^N}}}{\mathbf{E}}\left[e^{\lambda\mathcal{L}_{P_\mathsf{K}}^{\ell}(\mathfrak{f}_{\scalebox{.6}{$\uptheta$}})-\mathcal{L}_{\mathbb{D}'}^{\ell}(\mathfrak{f}_{\scalebox{.6}{$\uptheta$}})}\right]\right)}{\lambda}$,\\
$\frac{R}{K}=R_{\lambda'=\lambda K}= \frac{\mathbf{KL}(Q_{\text{\scalebox{.8}{\tiny{REV}}}}||P_{\text{\scalebox{.8}{\tiny{REV}}}})+\log\left(\frac{1}{\uptau}\right)+\log\left(\underset{\substack{{\hat\theta \sim P_{\scalebox{.8}{\text{\tiny{PAC}}}}},{\mathbbmsl{D}'\sim (P_\mathsf{K})^N}}}{\mathbf{E}}\left[e^{\lambda'\mathsf{K}^{-1}\mathcal{L}_{P_\mathsf{K}}^{\ell}(\mathfrak{f}_{\scalebox{.6}{$\uptheta$}})-\mathcal{L}_{\mathbb{D}'}^{\ell}(\mathfrak{f}_{\scalebox{.6}{$\uptheta$}})}\right]\right)}{\lambda'}$.\\
\\
\\
Furthermore, Donsker-Varadhan's change of measure states that, for and measurable function $\Phi:\mathcal{F}\to \mathbb R$ we have
\begin{flalign*}
{\mathbf{E}}\left[\Phi(\mathfrak{f}_{\uptheta})\right]\leq \mathbf{KL}(Q_{\text{\tiny{REV}}}||P_{\text{\tiny{REV}}}) + \log\left(\underset{{\hat\theta\sim Q_{\text{\scalebox{.8}{\tiny{REV}}}}}}{\mathbf{E}}\left[e^{\Phi(\mathfrak{f}_{\uptheta})}\right]\right).
\end{flalign*}
Thus, with $\Phi(\mathfrak{f}_{\uptheta})\stackrel{\mathrm{def}}=\lambda\left(\underset{{\hat\theta\sim Q_{\text{\scalebox{.8}{\tiny{REV}}}}}}{\mathbf{E}}\left[\mathcal{L}_{ P_\mathsf{K}}^{\ell_{\text{nll}}}(\mathfrak{f}_{\uptheta})\right]-\underset{{\hat\theta\sim Q_{\text{\scalebox{.8}{\tiny{REV}}}}}}{\mathbf{E}}\left[\mathcal{L}_{\mathbb{D}}^{\ell_{\text{nll}}}(\mathfrak{f}_{\uptheta})\right]\right)$, we obtain for all $Q_{\text{\scalebox{.8}{\tiny{PAC}}}}$ in $\mathcal{G}_{\hat{\Theta}}$:
\begin{flalign*}
&\lambda\left(\underset{{\hat\theta\sim Q_{\text{\scalebox{.8}{\tiny{REV}}}}}}{\mathbf{E}}\left[\mathcal{L}_{ P_\mathsf{K}}^{\ell_{\text{nll}}}(\mathfrak{f}_{\uptheta})\right]-\underset{{\hat\theta\sim Q_{\text{\scalebox{.8}{\tiny{REV}}}}}}{\mathbf{E}}\left[\mathcal{L}_{\mathbb{D}}^{\ell_{\text{nll}}}(\mathfrak{f}_{\uptheta})\right]\right)= \underset{{\hat\theta\sim Q_{\text{\scalebox{.8}{\tiny{REV}}}}}}{\mathbf{E}}\left[\lambda\left(\mathcal{L}_{P_\mathsf{K}}^{\ell_{\text{nll}}}(\mathfrak{f}_{\uptheta})-\mathcal{L}_{\mathbb{D}}^{\ell_{\text{nll}}}(\mathfrak{f}_{\uptheta})\right)\right]\\
&\leq \mathbf{KL}(Q_{\text{\tiny{REV}}}||P_{\text{\tiny{REV}}})+\log\left(\underset{\substack{{\hat\theta \sim P_{\scalebox{.8}{\text{\tiny{PAC}}}}}}}{\mathbf{E}}\left[e^{\lambda(\mathcal{L}_{P_\mathsf{K}}^{\ell}(\mathfrak{f}_{\scalebox{.6}{$\uptheta$}})-\mathcal{L}_{\mathbb{D}'}^{\ell}(\mathfrak{f}_{\scalebox{.6}{$\uptheta$}}))}\right]\right).
\end{flalign*}
Now, we apply Markov's inequality on the random variable $\zeta_{P_{\scalebox{.8}{\text{\tiny{PAC}}}}}(\mathbb{D}')\stackrel{\mathrm{def}}= \underset{\substack{{\hat\theta \sim P_{\scalebox{.8}{\text{\tiny{PAC}}}}}}}{\mathbf{E}}\left[e^{\lambda(\mathcal{L}_{P_\mathsf{K}}^{\ell}(\mathfrak{f}_{\scalebox{.6}{$\uptheta$}})-\mathcal{L}_{\mathbb{D}'}^{\ell}(\mathfrak{f}_{\scalebox{.6}{$\uptheta$}}))}\right]$ and moreover define $\zeta_{P_{\scalebox{.8}{\text{\tiny{PAC}}}}}(\mathbb{D}')_{\mathsf{K}}=\underset{\substack{{\hat\theta \sim P_{\scalebox{.8}{\text{\tiny{PAC}}}}}}}{\mathbf{E}}\left[e^{\frac{\lambda}{_{\mathsf{K}}}(\mathcal{L}_{P_\mathsf{K}}^{\ell}(\mathfrak{f}_{\scalebox{.6}{$\uptheta$}})-\mathcal{L}_{\mathbb{D}'}^{\ell}(\mathfrak{f}_{\scalebox{.6}{$\uptheta$}}))}\right]$:
\begin{align*}
&(P_{\mathsf{K}})^N\left(\zeta_{P_{\scalebox{.8}{\text{\tiny{PAC}}}}}(\mathbb{D}')\leq\frac{1}{\uptau}\underset{{\mathbbmsl{D}'\sim (P_\mathsf{K})^N}}{\mathbf{E}}\left[\zeta_{P_{\scalebox{.8}{\text{\tiny{PAC}}}}}(\mathbb{D}')_{\mathsf{K}}\right]\leq \frac{1}{\uptau} \mathcal{L}\left(\lambda\mathsf{K}^{-1}\right)\right)\\
&\geq 1-\uptau\frac{\underset{{\mathbbmsl{D}'\sim (P_\mathsf{K})^N}}{\mathbf{E}}\left[\zeta_{P_{\scalebox{.8}{\text{\tiny{PAC}}}}}(\mathbb{D}')\right]}{\mathcal{L}\left(\lambda\mathsf{K}^{-1}\right)}\\
&\geq 1-\uptau\frac{\mathcal{L}\left(\lambda\right)}{\mathcal{L}\left(\lambda\mathsf{K}^{-1}\right)},
\end{align*}
where $\mathcal{L}\left(\lambda\right)$ can be found in Appendix~\ref{pacbayesian},  Equation~(\ref{L2}).\\
This implies that with probability at least $1-\uptau\frac{\mathcal{L}\left(\lambda\right)}{\mathcal{L}\left(\lambda\mathsf{K}^{-1}\right)}$: 
\begin{flalign*}
\underset{{\hat\theta\sim Q_{\text{\scalebox{.8}{\tiny{PAC}}}}}}{\mathbf{E}}\left[\mathcal{L}_{P_\mathsf{K}}^{\ell}(\mathfrak{f}_{\uptheta})\right]-\underset{{\hat\theta\sim Q_{\text{\scalebox{.8}{\tiny{PAC}}}}}}{\mathbf{E}}\left[\mathcal{L}_{\mathbb{D}}^{\ell}(\mathfrak{f}_{\uptheta})\right]\leq \frac{1}{\lambda}\left(\mathbf{KL}(Q_{\text{\scalebox{.8}{\tiny{PAC}}}}||P_{\text{\scalebox{.8}{\tiny{PAC}}}})+\log\left(\frac{1}{\uptau}\right)+\underset{{\mathbbmsl{D}'\sim (P_\mathsf{K})^N}}{\mathbf{E}}\left[\zeta_{P_{\scalebox{.8}{\text{\tiny{PAC}}}}}(\mathbb{D}')_{\mathsf{K}}\right]\right),
\end{flalign*}
This statement is new, but main parts are derived from~\cite[Appendix A.2]{germain2016pac}. The extension follows for $\mathfrak{f}_{\uptheta}=\mathfrak{f}_{{\mathsf{k},\uptheta}}$ 

\end{document}